\documentclass[3p,twocolumn]{elsarticle}
\usepackage{lineno,hyperref}
\modulolinenumbers[5]

\usepackage[numbers]{natbib}
\usepackage[T1]{fontenc} 
\usepackage{fourier} 
\usepackage[english]{babel} 
\usepackage{amsmath,amsfonts,amsthm} 
\usepackage{caption}
\usepackage{subcaption}
\usepackage{graphicx}
\usepackage{wrapfig}
\usepackage{pdfpages} 
\usepackage[super]{nth}
\usepackage{calc}
\usepackage{float}
\usepackage{blindtext} 
\usepackage{hyperref}  
\usepackage{cleveref}
\usepackage{listings}   
\usepackage{booktabs}
\usepackage{enumitem}

\usepackage{fancyhdr} 
\usepackage[none]{hyphenat}

\usepackage{setspace}
\usepackage{titlesec}

\renewcommand{\paragraph}[1]{\noindent\textbf{#1}\ \ \ }

\date{\normalsize\today} 

\begin{document}

\begin{frontmatter}

\title{Automatic Recognition of Landmarks on Digital Dental Models}


\author[UoD-maths]{Br\'enainn Woodsend} 
\ead{bwoodsend@gmail.com}

\author[UoD-dental]{Eirini Koufoudaki}
\ead{e.koufoudaki@dundee.ac.uk}


\author[UoD-dental]{Peter A. Mossey}
\ead{p.a.mossey@dundee.ac.uk}

\author[UoD-maths]{Ping Lin\corref{cor1}}
\cortext[cor1]{Corresponding author}
\ead{p.lin@dundee.ac.uk}

\address[UoD-maths]{School of Science and Engineering, University of Dundee, Nethergate, Dundee DD1 4HN, United Kingdom}

\address[UoD-dental]{School of Dentistry, University of Dundee, Nethergate, Dundee DD1 4HN, United Kingdom}


\begin{abstract}

Fundamental to improving  Dental and Orthodontic treatments is the ability to quantitatively assess and cross-compare their outcomes. Such assessments require calculating distances and angles from 3D coordinates of dental landmarks. The costly and repetitive task of hand-labelling dental models impedes studies requiring large sample size to penetrate statistical noise.

We have developed techniques and software implementing these techniques to map out automatically, 3D dental scans. This process is divided into consecutive steps -- determining a model's orientation, separating and identifying the individual tooth and finding landmarks on each tooth -- described in this paper. Examples to demonstrate techniques and the software and discussions on remaining issues are provided as well. The software is originally designed to automate Modified Huddard Bodemham (MHB) landmarking for assessing cleft lip/palate patients. Currently only MHB landmarks are supported, but is extendable to any predetermined landmarks. 


This software, coupled with intra-oral scanning innovation, should supersede the arduous and error prone \textit{plaster model and calipers} approach to Dental research and provide a stepping-stone towards automation of routine clinical assessments such as "index of orthodontic treatment need" (IOTN).

\end{abstract}

\begin{graphicalabstract}
\centering
\includegraphics[width=.8\textwidth]{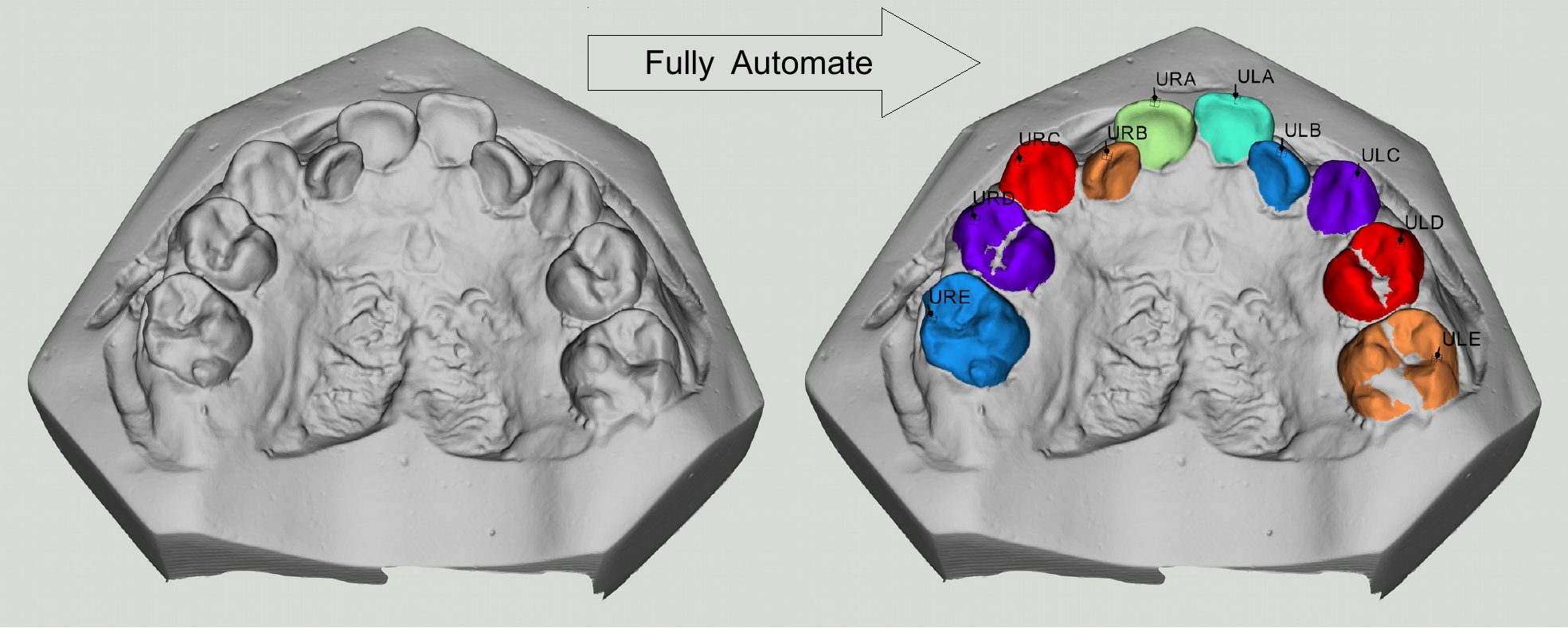}
\end{graphicalabstract}


\end{frontmatter}

\section{Introduction}

Three dimensional digital imaging has emerged as a new tool in clinical practice, and provides opportunities for improving research in multiple directions. In dentistry specifically, 3D analysis of jaws and dentition for treatment planning and treatment outcome assessment has already been gold standard (significantly before the digital era) in daily practice, especially in disciplines like orthodontics. 
That has been achieved through dental impressions followed by construction of three-dimensional plaster dental models. These models are then analysed by identifying on them landmarks and measurements of specific parameters. That manual process has always been problematic, as it is time consuming and is subject to both random and systematic errors. 
Recently intra-oral scanners have been developed that can deliver high accuracy digital models of single teeth and full dental arches. These digital models don't require storage, can be shared without needing to be shipped, and by annotating landmarks, lend themselves to digital analysis. \cite{Forsyth1996-lb,Rudolph1998-ac,Y_J_Chen_S_Kuang_C_Dds_H_Fu_C_Dds_and_K_Chee2000-zf}.
In this paper, we explore techniques to automate the finding of these landmark features on a digital 3D model of sets of teeth and the landmark 3D coordinates. 

Accurate automated dental landmark identification would be a great tool both to researchers of dental science, and in routine treatment planning and assessment in clinical dentistry. Tooth measurements have always been regarded as time-consuming \cite{Knyaz2016-gd}, and so automatic landmarking would save time for dentists, and opens up the possibility of studies on large numbers of teeth sets.

\begin{figure}[ht]
    \centering
    \captionsetup{margin=.04\textwidth}
	\includegraphics[width=0.95\linewidth]{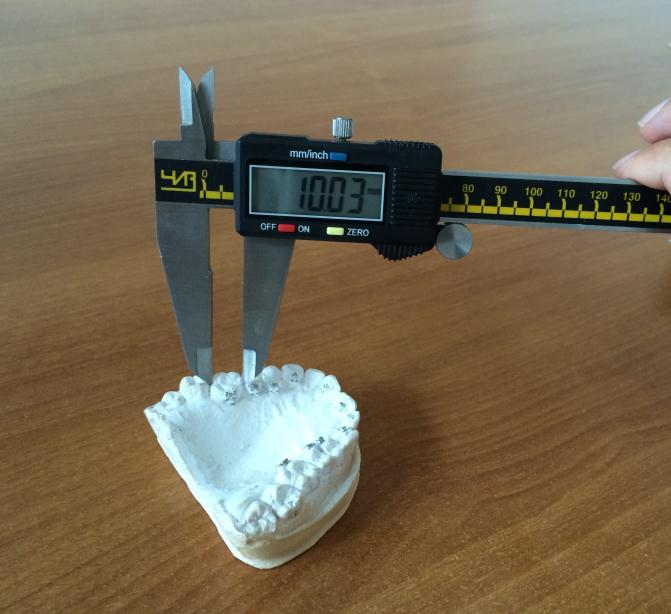}
    \caption{Manual measurements on plaster models}
    \label{fig:manual-measuring}
\end{figure}

\section{Background}

Traditionally dental impressions are a negative imprint taken from the patients mouth using an impression tray and a setting paste. By many patients this is considered a very uncomfortable procedure. Then the negative imprints are poured with plaster, which then sets creating 3D plaster copies of dental arches. These arch copies are then used to measure and estimate tooth position in relation to neighbouring and opposing teeth (Figure~\ref{fig:manual-measuring}). 
The development of digital impressions have been a game changer, as it can develop dental arch imprints of high accuracy very quickly and with greater comfort for the patient. The digital data can either be collected directly from the patient or by creating a plaster model first that can then be scanned. Of course the first option is much more efficient in multiple ways, including time, materials, storage management and human resources. 

The 3D scans produced by the intra-oral scanners are in STL format -- a standard open source format for 3D models. They describe only the surface of a model and therefore require wholly different analysis techniques to voxel based \textbf{volumetric} scans such as the output of a micro-CT or 2D images such as X-ray. 

The original objective of this study is to automate the modified Huddart and Bodenham system (MHB) \cite{Mossey2003-lg} for assessing treatment of cleft lip/palate. This system has been shown to be far more objective and reproducible than its predecessors \cite{Gray2005-zs}, but has had limited uptake by clinicians due the the considerably extended time it takes. These traits make it a prime target to be converted into an automated software. Ma~et~al~\cite{Ma2017-em} devised a semi-automatic system, automatically calculating an MHB score based on manually selected coordinates of key landmarks. These particular landmarks are the midpoints of the incisors, the tips of the canines and the outer cusps (bumps) on the molars. Manual landmark placing is time consuming, requires expertise and is prone to human error. We have created a software which is set to identify dental landmarks in accordance with the MHB scoring system (although the software can be adjusted to work with many systems). The aim of the application is to increase efficiency and automation of the scoring of dental surgical outcomes, encouraging a more efficient workforce in global dental care. By moving from traditional plaster ``hard-copy'' models to 3D digital models, the global burden of care will be reduced. In addition, the reliability and reprehensibility of dental model scoring will improve by reducing human error and increasing the accuracy of measurements.

\subsection{Similar prior arts}
There have been several studies that have worked towards similar or overlapping objectives as those of this work. Some of the most interesting and partially relatable ones to our development are summarised here. 

One of the oldest systems, developed by Kumar~et~al~ \cite{Automatic-Feature-Identification-in-Dental-Meshes}, aims for fully autonomous mapping out of a dentition in a slightly different methodology than the one followed by us. An orientation step is not part of the procedure as it is  assumed that all models are oriented in the same way. Then they introduce a \textit{watershed} method to partition the teeth. This method is analogous of flooding a mesh with water until small lakes are formed -- these lakes are the teeth -- with the catch that the \textit{height} of a mesh vertex is in-fact defined based on mixture of surrounding curvature as well as its regular geometric height. Teeth are then identified using curvature of cross sections. The watershed method has persisted, being reused in many more modern works. 

Considerable work has been put into tooth segmentation, often with the \textit{fully automatic} constraint relaxed to \textit{semi-automatic}. This is largely driven by forensic scientists who wish to identify postmortems where dentitions were only partially recovered. Kronfeld~et~all~\cite{snake} have designed a \textit{snake} algorithm to \textit{walk} along the edges of teeth to partition them, with some safeguards to help jump gaps in tooth outlines. Zou~et~al~\cite{harmonic-field-ZOU2015132} have developed a tooth partitioning system based on user-supplied tooth labels. It turns the mesh into a graph network (dubbed \textit{a harmonic field}) with each vertex a node and each edge an arc. Each unlabeled node has an unknown potential and each labelled node has a fixed potential dependent on its label. Each arc has a weight derived from curvature, and flow through it proportional to the difference in potentials of its two vertices and inversely proportional to its weight. The whole system of vertex potentials/arc flows is solved with the constraint: the net flow of each non-labelled vertex must be zero. A vertex is part of a tooth if its potential exceeds some halfway threshold. Each tooth must be separated one at a time, but if you solve the linear system using sparse LU factorisation, the potentials of each labelled vertex can be altered with negligible extra computation i.e. the processor-intensive part is done only once. This method truly tackles the issue of \textit{shabby} casts where the tooth edges are poorly defined or have gaps.

Lastly, Kalogerakis~et~al~\cite{Kalogerakis:2010:labelMeshes} focuses on simultaneous segmentation and labelling of arbitrary objects (not dentally related) using machine learning. They manage it with great success using training sets of less than 10 models per object type.

\section{Materials/Methods}
\label{sec:automated-landmarks}

The software was developed in the programming language Python. We had available 239 dental models of different types (listed in table~\ref{tab:models-counts}). The models are stored in STL format. 

An STL file describes the only surface of a 3D object, making it hollow. It contain no colour or texture information. Once read from file, it is typically referred to as a \emph{mesh} - an unordered list of triangles with each triangle defined by the $(X, Y, Z)$ values of its three vertices. Since all triangles are represented separately, vertices that are corners of multiple triangles are duplicated. A typical first step on reading an STL, which our software adopts, is to find and enumerate all the unique vertices to make the connections between neighbouring triangles easier to find. Traditionally, STL files contain no scale information and their units are generally arbitrary, but in dental scans the units are always millimetres. 

\begin{table}[!ht]
\centering
\begin{tabular}{rlll}
\hline
Count   & Qualifiers   \\
\hline
 24 & Upper &         &              \\
 21 & Lower &         &              \\
 17 & Lower & IO      &              \\
 16 & Upper & IO      &              \\
 81 & Upper & Primary & Cleft        \\ 
 80 & Lower & Primary & Cleft        \\
\hline
\end{tabular}
\caption{Types and counts of models used in this study}
\label{tab:models-counts}
\end{table}

\begin{figure*}[p]
    \vspace*{-2cm}
    \makebox[\linewidth]{
        \includegraphics[scale=.8]{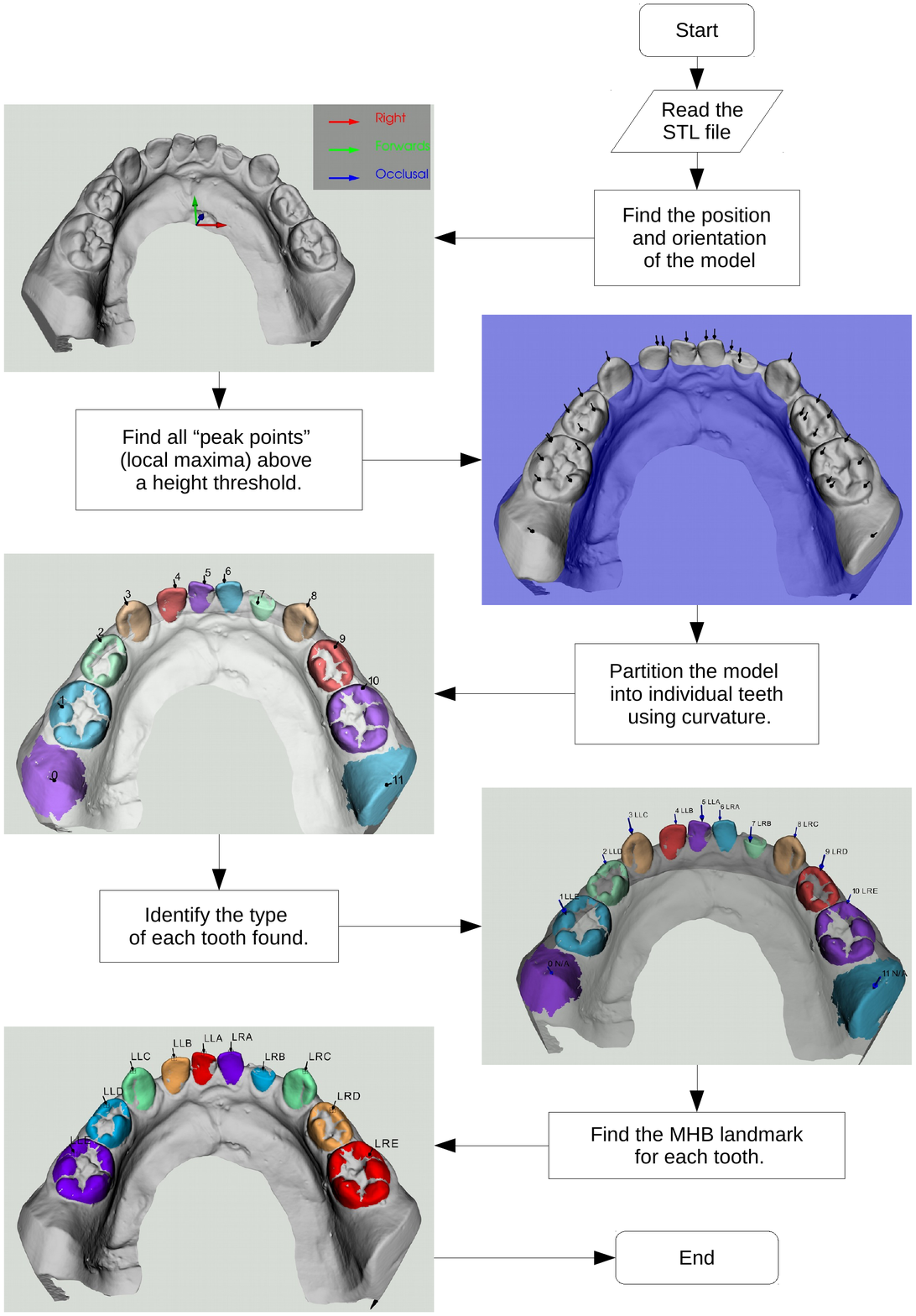}
    }
    \caption{Flowchart summary of the automatic landmark finding process}
    \label{fig:main-methods}
\end{figure*}

The data analysis for the successful automated landmark identification, is a multi step process approach. The steps undertaken by the software are summarised by the flow chart in figure~\ref{fig:main-methods}.

\subsection{Find Orientation}
\label{sec:orientation}

\newcommand{\unitvec}[1]{\hat{\mathbf{#1}}}
\newcommand{\evector}[1]{\unitvec{e}_{\text{#1}}}

Before any kind of analysis can be done, the model's orientation needs to be established. As different scanners use wholly different orientation conventions and initial the orientations aren't accurate to begin with, it is important that the system is able to find and internally normalise the orientation of models before they can be further analysed.

Orientation is expressed as a series of \emph{unit-vectors} representing \emph{left}, \emph{right}, \emph{backwards}, \emph{forwards}, \emph{down}, \emph{up} and \emph{occlusal}. These unit-vectors should be derived from the model and any inspection of a vertex should use these vectors so that the effect of the model's position is negated. i.e. To get the height of a vertex take the inner product of that vertex with the \emph{up} vector instead of examining raw Z values (or Y values on a TRIOS-scanned model). The model itself must not move or the relative positioning between a patient's upper and lower jaws will be lost. 

Section~\ref{sec:height-threshold} requires that the vertical directions are particularly accurate. The horizontal directions are less crucial. To generate the required unit-vectors a PCA (principle component analysis -- see section~\ref{sec:pca}) based method was designed.

\subsubsection{Principal Component Analysis}
\label{sec:pca}

PCA looks at covariance (spread from the centre of mass\footnote{The centre of mass being the middle of the model, or more precisely, the mean of all of its vertices.}) in all directions. It returns a set of axes (three perpendicular unit-vectors) ordered from most to least covariance. PCA may be used to find the directions of the longest, middle-length, and shortest dimensions of an object. PCA calculations are available in the appendix~\ref{sec:calculate-pca}.

A dental model is wider (left/right) than it is long (forwards/backwards) and longer than it is tall. So the output of PCA on a dental model should yield \emph{left}/\emph{right} as the first (largest covariance) unit-vector, \emph{forwards}/\emph{backwards} as the second and  \emph{up}/\emph{down} as the third. 

PCA’s advantages are: (i) it is one of most popular statistical techniques and thus the function of implementing 
it is directly available in any common programming environment; (ii)  the model can be positioned and orientated anywhere and PCA will track it indifferently.

\subsubsection{Signs of the unit-vectors}
\label{unit-vector-signs}

\begin{figure}[ht]
    \centering
    \captionsetup{margin=.025\linewidth}
    \includegraphics[width=.95\linewidth]{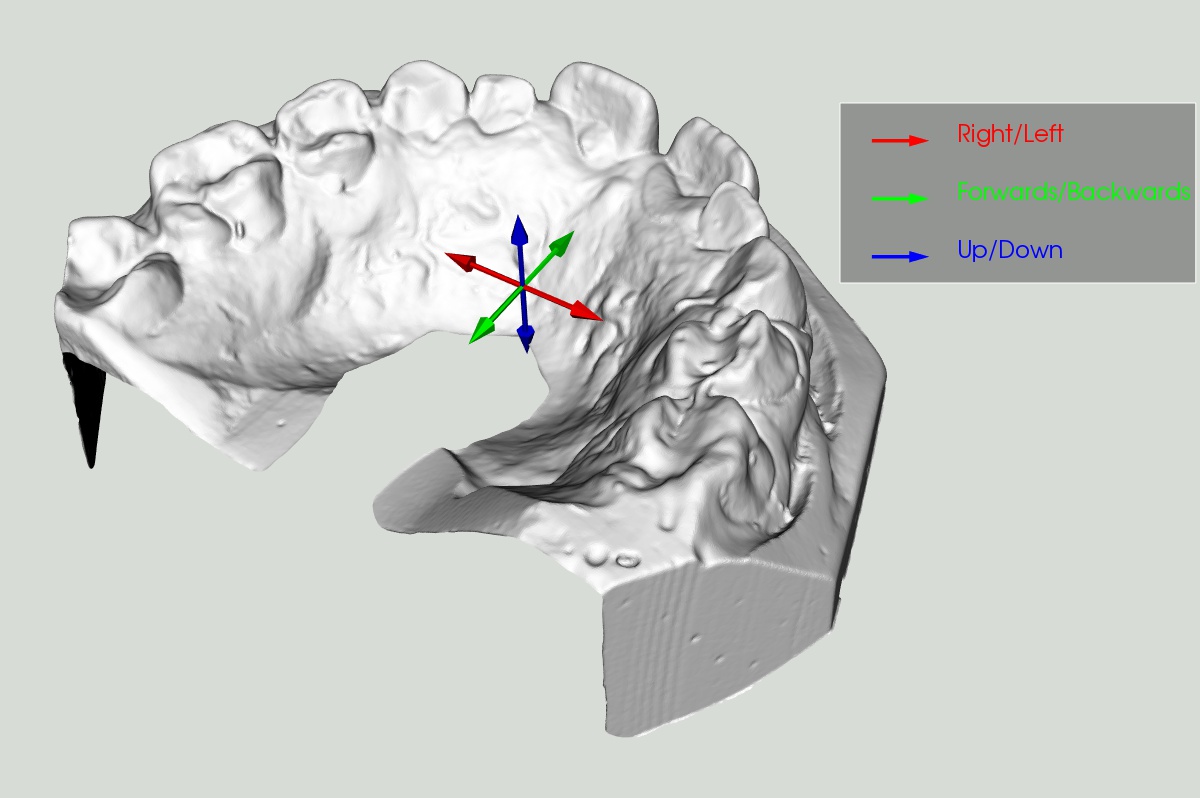}
    \caption{The unit-vectors from PCA. The red arrows represent \emph{left} and \emph{right}, but PCA doesn't specify which is \emph{left} and which is \emph{right}. Likewise with \emph{forwards}/\emph{backwards} (green) and \emph{up}/\emph{down} (blue).}
    \label{fig:odom-unsigned}
\end{figure}

The sign of an eigenvector is arbitrary and consequently so are those of the unit-vectors found using the above (see figure~\ref{fig:odom-unsigned}). The signs have to be checked by other means and the vectors reversed if they are wrong. The following checks were adopted (in the following order).

\paragraph{1. Vertical \emph{up}/\emph{down}/\emph{occlusal}}
As the occlusal surface is the most detailed, the density of triangles (and their corresponding unit-normals) there is highest. A mean of all the model’s unit normals should therefore approximately point \emph{occlusally} and PCA’s \emph{occlusal} can be sign matched to this approximate \emph{occlusal}.

\begin{figure}[ht]
    \centering
    \captionsetup{margin=.025\linewidth}
    \includegraphics[width=.95\linewidth]{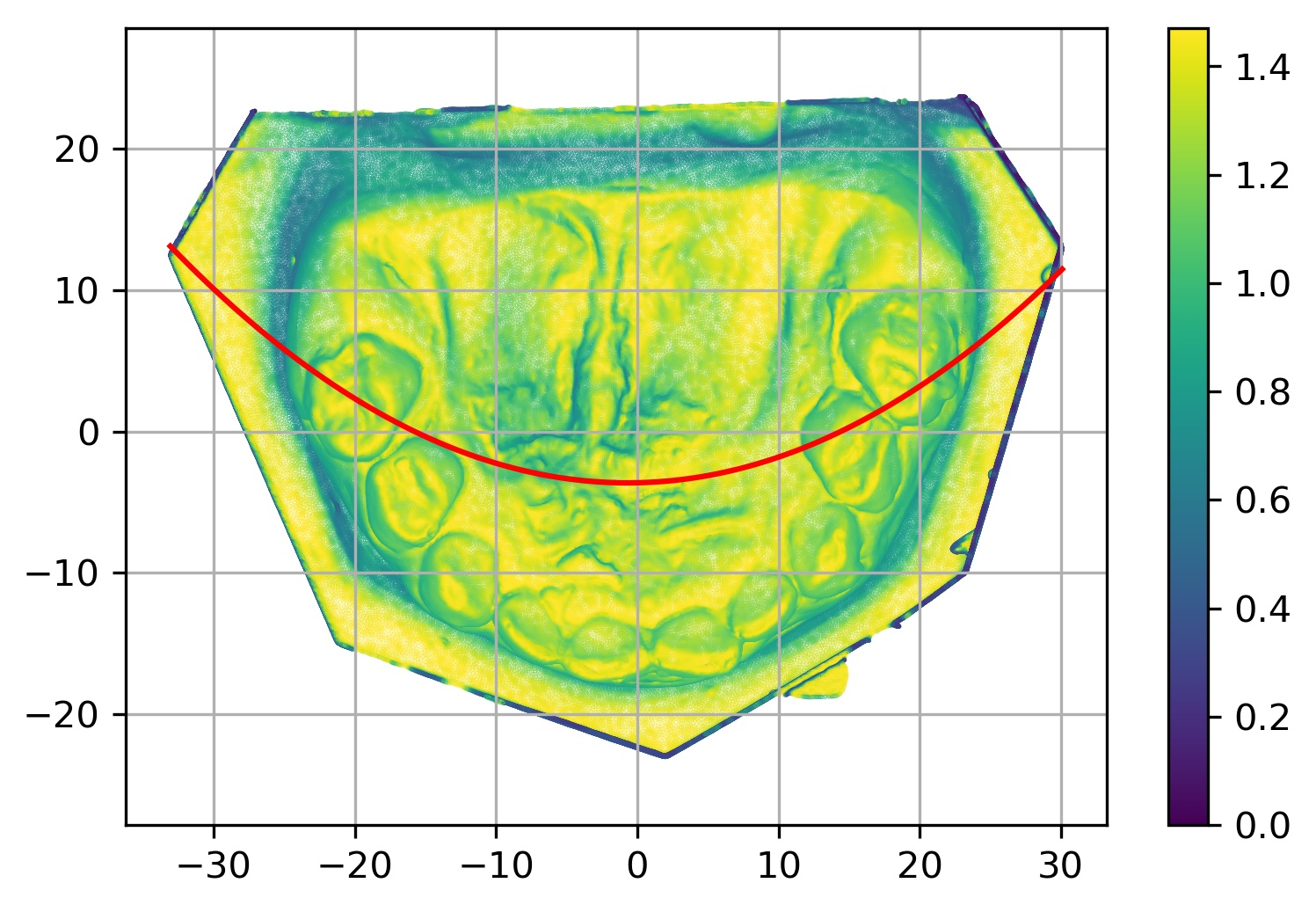}
    \caption{Quadratic weighted fit (red line) used to test the sign of the \emph{forwards} unit vector. The line's being $\cup$ shaped tells us that \emph{forwards} is actually \emph{backwards} and needs reversing.}
    \label{fig:odom-y}
\end{figure}

\paragraph{2. Forwards/Backwards}
The software extracts the horizontal components of the mesh's vertices and approximately captures the \emph{arch} shape of the jaw. All the points on the mesh are reverse weighted by their unit normals’ agreement with \emph{occlusal}. This is supposed to emphasise the labial and lingual surfaces. A weighted quadratic curve is fitted to the horizontal components. If the \emph{forwards} vector’s sign is correct then the quadratic should be $\cap$ shaped with a negative $x^2$ coefficient. If it is $\cup$ shaped then the \emph{forwards} vector needs reversing.

\paragraph{3. Left/Right}
With the signs of the other two axes known this can just be determined so as not to mirror the model. Treating the unit-vectors as column vectors, mirroring can be checked for using:

\begin{equation}
    \text{det} \:\:
    \begin{bmatrix} \evector{right} & \evector{forwards} & \evector{up} \end{bmatrix}
     \:\: = \:\: \left\lbrace \begin{aligned}
     1,  \quad & \text{Non-mirroring}\\
    -1 \quad & \text{Mirroring}\\
    \end{aligned} \right.
\end{equation}

If it mirrors then reverse $\evector{right}$.

\subsubsection{Fine-Tuning the Vertical Axis}

\begin{figure}[ht]
    \centering
    \begin{subfigure}{0.45\textwidth}
        \captionsetup{margin=.05\textwidth}
        \includegraphics[width=.95\textwidth]{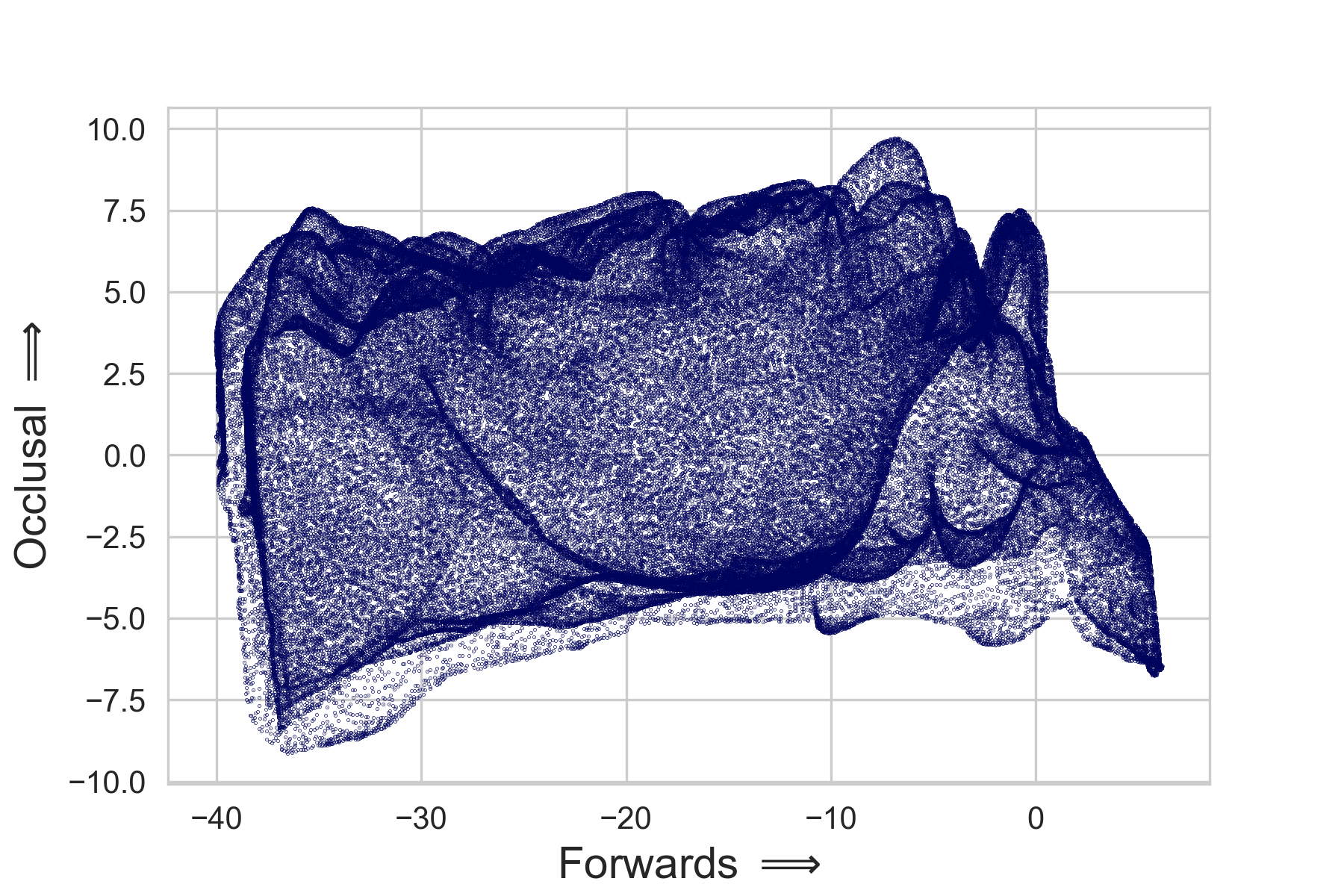}
        \caption{Side-on flattened screenshot of all the points \newline on a model. This is a lower jaw with the patient looking to the right.}
        \label{fig:odom-side1}
    \end{subfigure}
    \begin{subfigure}{0.45\textwidth}
        \captionsetup{margin=.05\textwidth}
        \includegraphics[width=.95\textwidth]{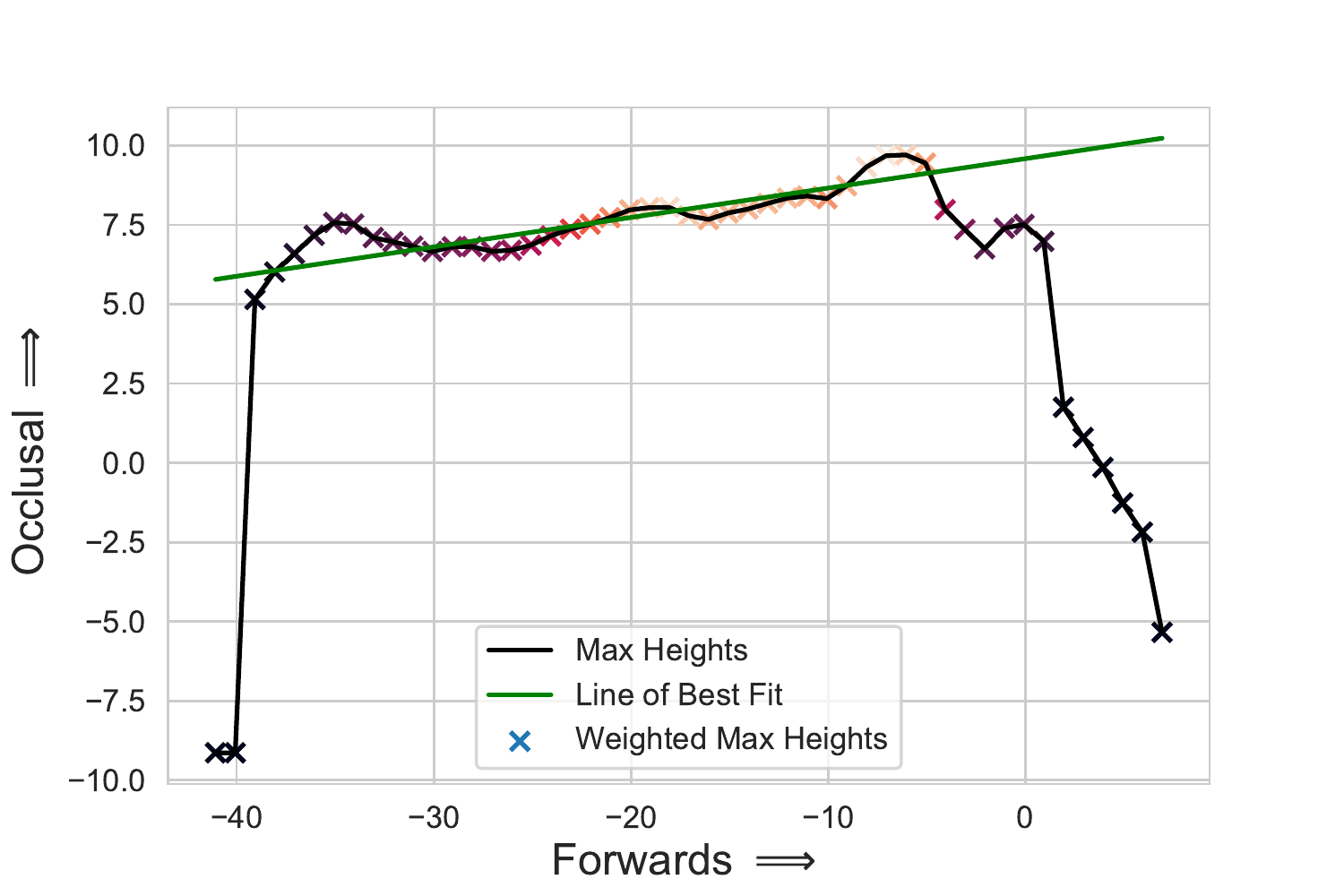}
        \caption{The top-most outline (black line) of the same side on view, the weighted points along that outline ('x' shaped markers) and the line fitted to those points (green line) which will represent \emph{horizontal}.}
        \label{fig:odom-side2}
    \end{subfigure}
    \label{fig:oom-side}
    \caption{Method to fine-tune the vertical axis}
\end{figure}

Later steps (mostly section~\ref{sec:peak-points}) require particularly accurate vertical unit-vectors. Depending on the type of dental model (intra-oral scans and plaster models with rough bases), the vectors from PCA are typically inadequate. Using the PCA \emph{forwards} and \emph{occlusal} unit-vectors (after the above sign checking), more accurate unit-vectors can be found by fitting a straight line through the top-most outline (see figures~\ref{fig:odom-side1} and~\ref{fig:odom-side2}). To do this, divide the mesh horizontally into bins and find the highest point in each bin. Then fit to those highest points, aggressively weighting the centre-most (horizontally) and highest points so as not to be affected by dips due to missing teeth, or the drop-off at the front and back of the model.

\subsection{Find the Peak Points}
\label{sec:peak-points}

\begin{figure*}[!ht]
    \centering
    \begin{subfigure}{0.3\textwidth}
        \includegraphics[width=0.99\linewidth]{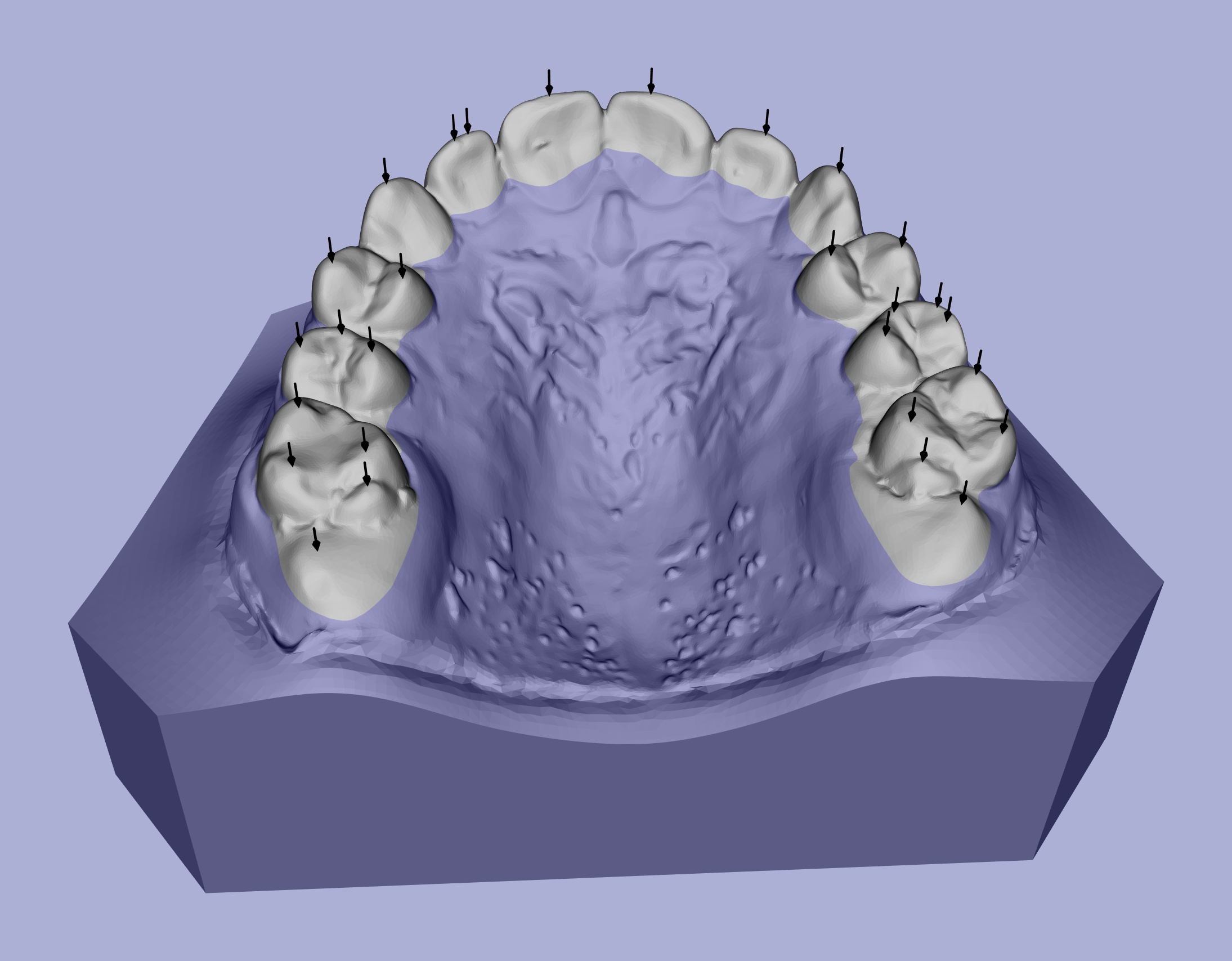}
        \caption{Adult Maxillary}
    \end{subfigure}
    \begin{subfigure}{0.3\textwidth}
        \includegraphics[width=0.99\linewidth]{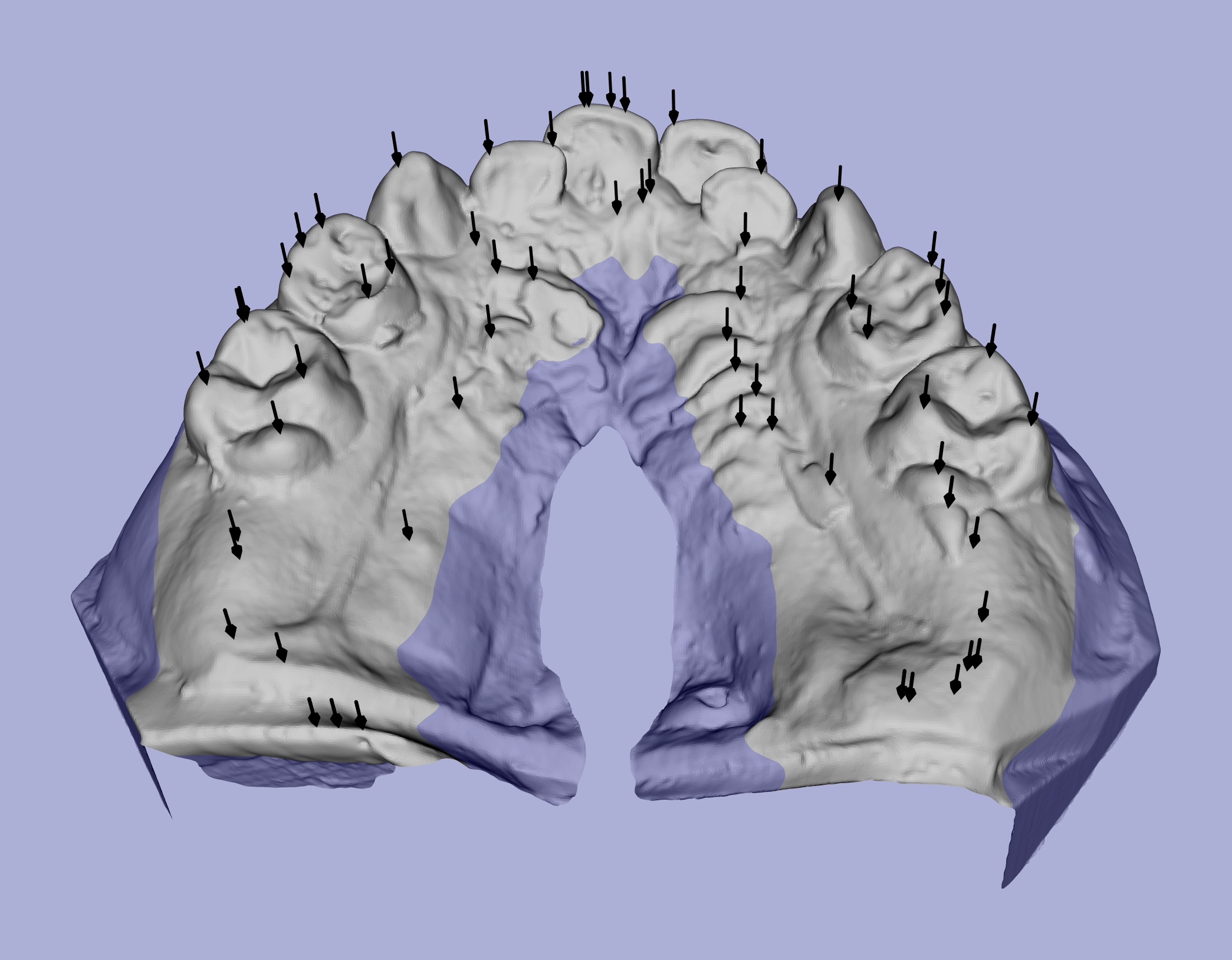}
        \caption{Primary (5 years old) Maxillary}
    \end{subfigure}
    \begin{subfigure}{0.3\textwidth}
        \includegraphics[width=0.99\linewidth]{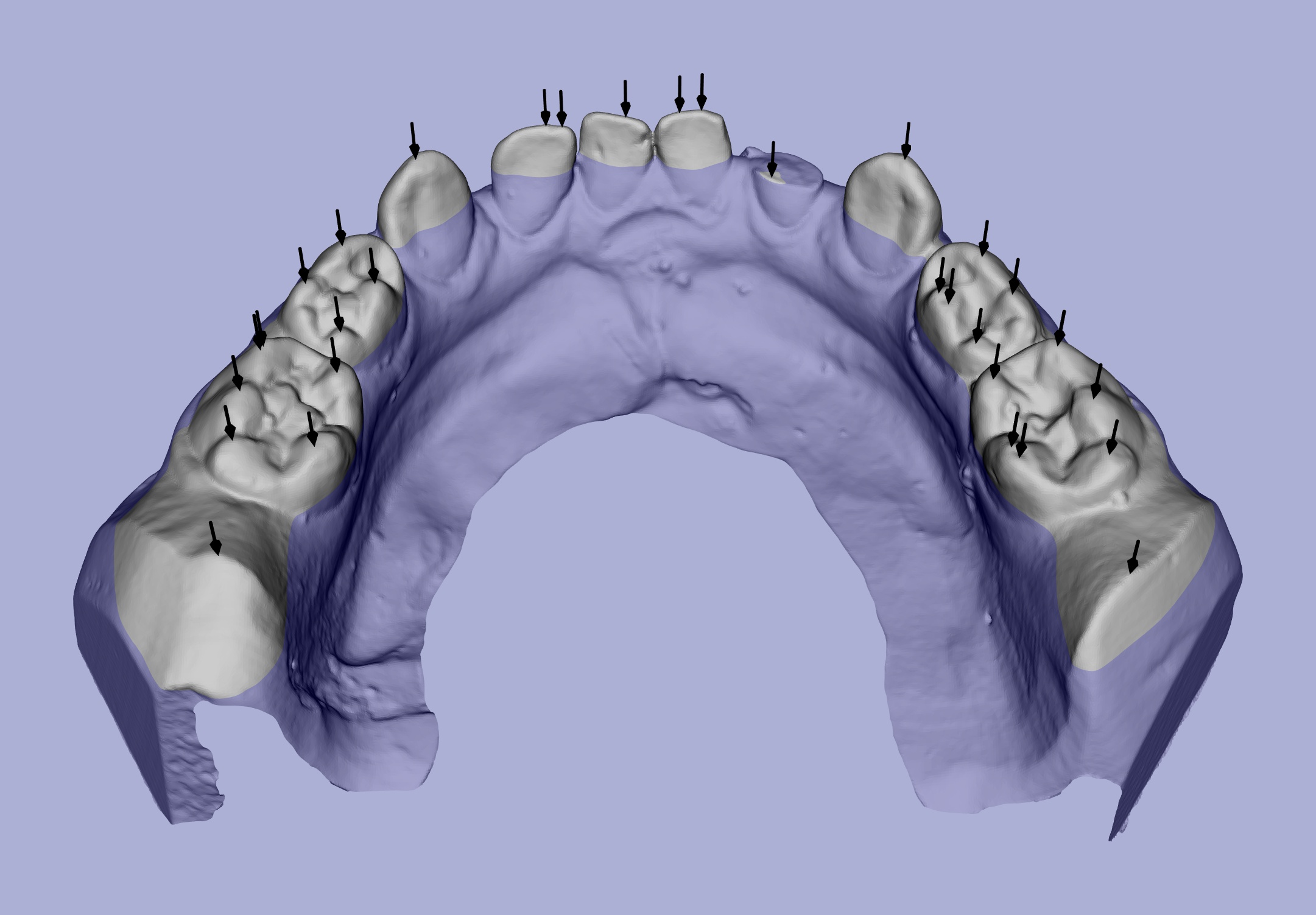}
        \caption{Primary (5 years old) Mandible}
    \end{subfigure}
    \newline
    \caption{The peak points (black arrows) and the height threshold (blue tint) found for three
    different models.}
    \label{fig:peak-points}
\end{figure*}

Most of the landmarks required for MHB scoring are either located on teeth tips or cusps. Mathematically these can be described as local maxima in the occlusal direction or any vertex that is higher than all its immediate neighbours. These points are referred throughout this article as \emph{peak points} or just \emph{peaks} (see figure~\ref{fig:peak-points}). No information about what feature each peak point represents is found at this stage. There will be many peak points that do not represent an actual landmark and must be cleared away in later steps. The only requirement of this step is to land at least one point on each tooth or cusp (for molars and premolars).

\hypertarget{sec:height-threshold}{
\subsubsection{Filtering Based on Height}\label{sec:height-threshold}}

To avoid finding large numbers of peaks on the gums, base (for plaster models) and roof of the mouth, the search area is reduced to only the top (occlusally) of the model so as to mostly include only teeth. Any points 6mm or more below the highest tip of the teeth are excluded from the search area. The height threshold is visualised with the transparent blue planes in figure~\ref{fig:peak-points}.

The height threshold chosen ($6mm$ below the highest peak) is comfortably low enough to include the tips or cusps of all teeth of interest -- even if they are chipped (although these may be mistakenly rejected in section~\ref{sec:tooth-assignment} -- the tooth assignment stage). Setting it very low increases the number of non-tooth features picked up. These features will be safely removed later (in section~\ref{sec:cleaning}) but at considerable computational expense.

\subsection{Partition into Individual
Teeth}
\label{sec:partition-into-individual-teeth}

\subsubsection{The General Idea}
\label{sec:the-general-idea}

\begin{figure*}[ht]
  \begin{subfigure}[b]{0.45\textwidth}
    \includegraphics[width=.95\linewidth]{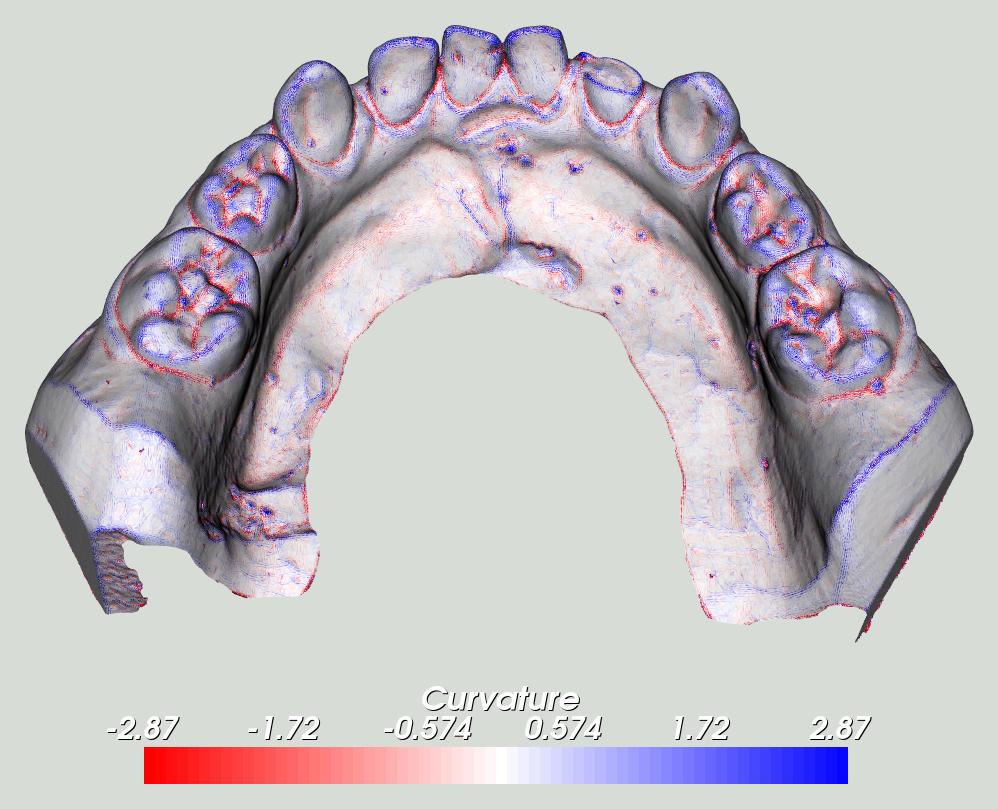}
    \caption{Front View}
    \label{fig:f1}
  \end{subfigure}
  \hfill
  \begin{subfigure}[b]{0.45\textwidth}
    \includegraphics[width=.95\linewidth]{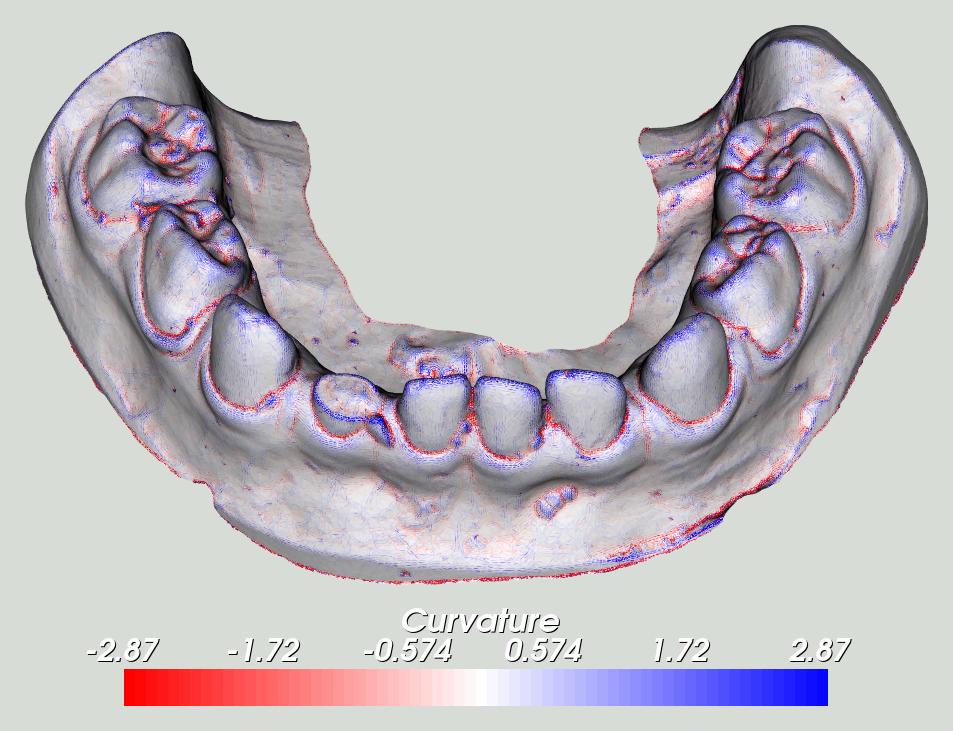}
    \caption{Rear View}
    \label{fig:f2}
  \end{subfigure}
  \caption{Signed Curvature -- Notice the ring of red along the tooth-gum join surrounding each tooth and between neighbouring teeth.}
  \label{fig:signed-curvature}
\end{figure*}

This step uses \emph{curvature} to find the boundaries of teeth.
Curvature is a quantitative measurement of how much a surface deviates
from being flat at a particular point. The exact definition of curvature can vary -- the one chosen is signed so that an outside corner (a bump, cap or tip) is positive and an inside corner (a slot, groove or crease) is negative.
The join between each tooth and the gum is a crease and therefore the curvature along the join is negative (see figure~\ref{fig:signed-curvature}).

 Starting at the top of a tooth and recursively including adjacent mesh triangles until an edge of significantly negative curvature is hit, one can find all triangles that are part of that tooth. This region covered will be referred to as the starting point's \emph{region}. Each peak point found in the last section is used as a starting point. Peaks on the same tooth will have regions that overlap. By testing for overlapping regions, duplicity of teeth is avoided. Any peaks that weren't on a tooth to begin with will not be bounded by the creases of tooth-gum joins and will therefore try to include most of the model if left unchecked. By imposing the rule \emph{stop if travelled more than a tooth's width away from the starting point}, and testing if that rule was actually used, non-tooth peaks and their corresponding regions can be identified. These non-tooth regions are labelled \emph{spilled} and are discarded.

Thus, each tooth should come out nicely partitioned without any duplicity. (No information about which tooth it is which is found here.) And all non-tooth features should remove themselves.

\paragraph{Complications} Whilst the above may seem promising it doesn't happen in practice. Below is listed some of the more prominent issues which must be solved.

\begin{itemize}
    \item Non-tooth features often don't spill. This happens mostly on the palatal rugae and on \textit{knobbly} plaster models.
    \item Each cusp of the molars and premolars will usually be separated and have to be re-grouped back together.
\end{itemize}

\begin{figure*}[p]
    \vspace*{-1cm}
    \makebox[\linewidth]{
        \includegraphics[width=\linewidth]{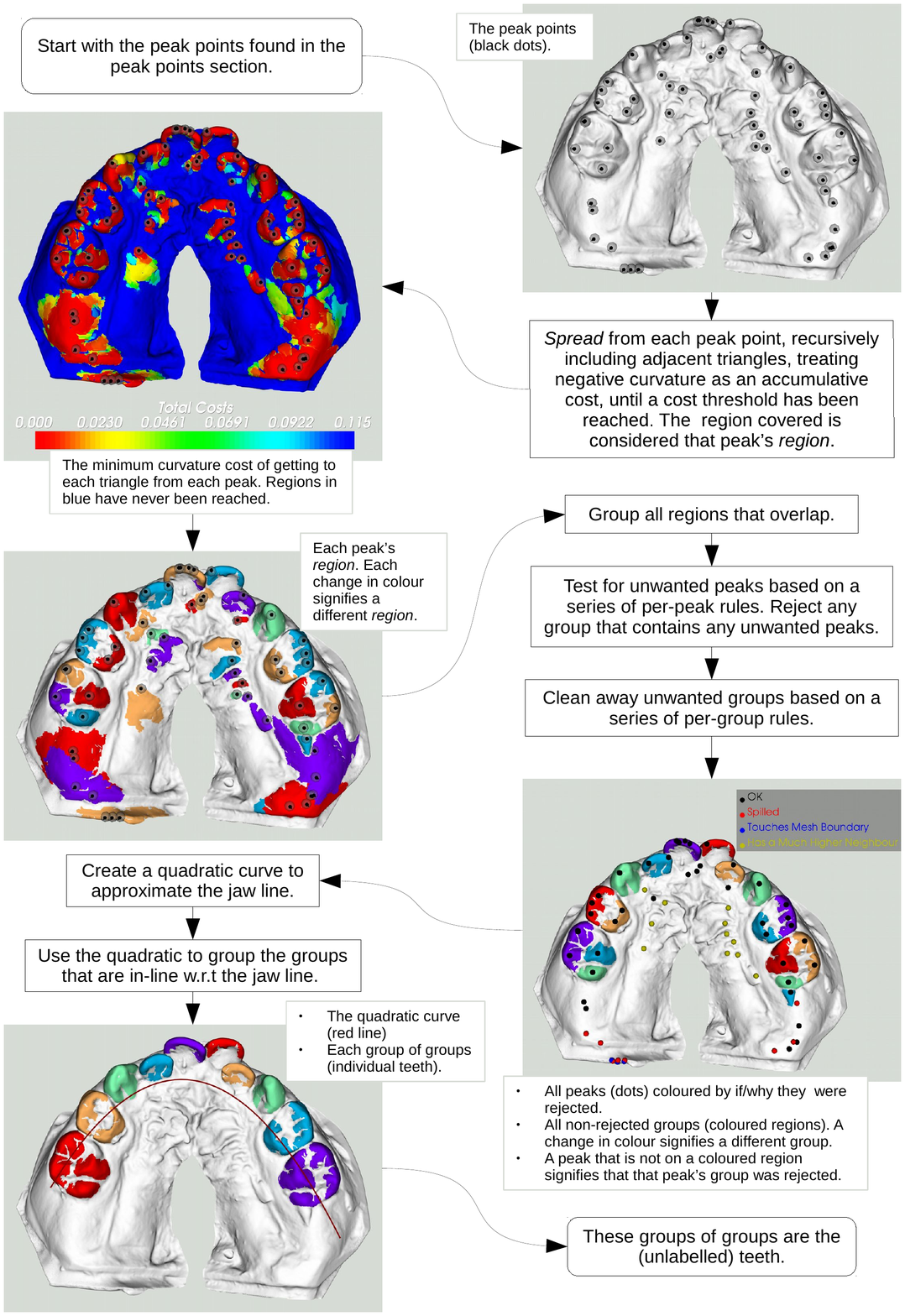}
    }
    \caption{Flowchart summary of the automatic landmark finding process}
    \label{fig:partition-flow}
\end{figure*}

The flowchart in figure~\ref{fig:partition-flow} outlines the full series of steps undertaken.

\subsubsection{Flood Filling From Each Peak-Point}

\paragraph{Calculating Curvature:} Curvature comes in many different forms. For continuously differentiable surfaces, it can be calculated using high derivatives of the surface geometry. But for a discretely defined surface, like that of an STL file, approximating high derivatives gives a very poor signal to noise ratio. Prior works have typically used Principal Curvatures which are defined per-vertex. Principal Curvatures are derived by either taking minimum \cite{4381713,single-tooth-modelling} or mean \cite{WU2014199,snake} of nearby edges' curvatures. 

Our system uses per-edge curvature rather than per-vertex curvature. Edge curvature was chosen as it is easier to calculate, involves slightly less approximating and it makes the after-analyses easier as each triangle always has exactly three edges and three neighbours (triangles in non-closed meshes will occasionally have less neighbours) whereas, with per-vertex curvature, each vertex can have any number of neighbours and edges, requiring awkward \textit{ragged-array} data structures.

Our form defines the curvature ($\mathbf{k}$) at an edge by comparing the outward unit normals ($\mathbf{\hat{n}}_0$ and $\mathbf{\hat{n}}_1$) of the two triangles on either side of that edge along with the distance $|\Delta \mathbf{x}|$ between the centres of each triangle, ($\Delta \mathbf{x}$ being the displacement from the centre of triangle 0 to the centre of triangle 1.)

\begin{equation*}
    \mathbf{k} = \frac{\unitvec{n}_0 \times \unitvec{n}_1} {|\Delta \mathbf{x}|}
\end{equation*}

This equation yields a vector, the direction of which is just the edge's direction (which we don't need). The magnitude quantifies how tight the corner of two mesh triangles sharing the edge is. This still leaves no information about the sign. The equation used to get the sign is displayed below:

\begin{equation*}
    \text{signed curvature} = - \text{sign}( \mathbf{\hat{n}}_0 \cdot \Delta \mathbf{x} ) \space |\mathbf{k}|
\end{equation*}

The output of the above equation is a signed curvature scalar, for each mesh triangle, for each of it's adjacent triangles (as shown in figure~\ref{fig:signed-curvature}).

\paragraph{A \emph{Cost Map}} is derived from the above signed curvature. Only creases (negative values) are relevant. And, being a \emph{cost}, all values should be positive.

\begin{equation*}
    \text{cost} = \text{max}( - \text{signed curvature}, \: 0)
\end{equation*}

Again, this cost map is per edge of each of the mesh's triangles. The cost of including a new adjacent triangle to the region is the cost of crossing the triangle edge separating that triangle from the region. Costs accumulate making a problem analogous to \emph{The Shortest Path Problem} from Graph Theory \cite{dijkstra1959note} with the exception that we are interested in the cheapest path to every mesh triangle rather than a single destination. In this analogy a mesh triangle is considered as a graph node and the weight or distance on the edge that connects two nodes (or triangles) is the cost defined above. Once the cumulative cost (i.e. the minimal cost solved from the shortest path problem) to every triangle is found it is compared to a \emph{cost threshold}. Any triangles with cumulative costs below the threshold are included as part of that peak's \textit{region}. The \textit{cost threshold} is solved for dynamically per model so as to maximize the total surface area classed as part of a tooth after the grouping and cleaning stages throughout the rest of this section.

\noindent Mathematically, the problem is formulated as follows. Let:

\begin{itemize}
    \item $T[i]$ be the total cost of reaching triangle $i$ (with $i\in[1,\text{number of triangles in the whole model}]$).
    \item $E[i,j] >=0$ be the curvature based edge cost of crossing the edge from triangle $i$ to triangle $j$, a neighbour of $i$.
    \item $T_{max} > 0$ be our curvature threshold.
\end{itemize}

\noindent Then solve the following for all elements of the vector $T$. 

\begin{itemize}
    \item $T[i] \:\: = \:\: 0 \quad$ 
        if any of triangle $i$'s vertices is the initial peak point.
    \item $T[i] \:\: = \:\: \min_{j}( T[j] + E[i, j]) \quad$
        if $\:\: \min_{j}( T[j] + E[i, j]) \:\: < \:\: T_{max}$.
    \item $T[i] \:\: = \:\: T_{max} \quad$ otherwise.
\end{itemize}

A triangle $i$ is then considered part of the region if $T[i] < T_{max}$.

Rather than truly solving the system as a linear algebra problem, which would be difficult due to the uses of $\min$, and slow due to the large number of triangles involved, a far more efficient algorithm was devised. This algorithm is very close to Dijkstra's algorithm \cite{dijkstra1959note} for solving the Shortest Path problem.

\begin{enumerate}
    \item First, initialise all $T[i]$s to $T_{max}$.
    \item Initialise an empty queue.
    \item For each triangle $i$ which contains the starting peak, set its $T[i]$ value to $0$ and add its three neighbouring triangles to the queue.
    \item Pop (choose and remove) an element $i$ from the queue.
    \item Evaluate $t = \min_{j}( T[j] + E[i, j])$ for that $i$ and its three neighbours $j$. It doesn't matter that one or more of the $T[j]$ may not have been processed yet. If $t < T[i]$ then set $T[i] = t$ and add the three neighbours $j$ to the queue to be (re)calculated later.
    \item If the queue is empty, terminate. Otherwise return to step 4.
\end{enumerate}

\subsubsection{Group Overlapping Regions}
\label{sec:group-overlapping}
Group all regions that overlap (i.e. have triangles in common). Indirect grouping is allowed, meaning that if regions $A$ and $B$ overlap, and regions $B$ and $C$ overlap, but $A$ and $C$ don't, then $A$, $B$ and $C$ should still form one group.

\subsubsection{Remove Non-Tooth Regions}
\label{sec:cleaning}

Non-tooth regions are tested for with a series of rules to remove unwanted features. It's ideal, but not imperative, that all non-tooth features are removed before proceeding with the tooth assignment section.

\paragraph{Per Peak or Per Peak's region Rules:}

A peak and the group it belongs to is rejected if it meets any of the following criteria:

\begin{itemize}

    \item It has \textit{spilled} (as defined in section~\ref{sec:the-general-idea}).

    \item It is in close horizontal proximity to a much higher peak. Or more precisely, if the ratio of the vertical displacement $\delta V$ and horizontal distance $\mod{\delta H}$ to any other peak is greater that $1.5$ then this peak is almost certainly gum.

\end{itemize}

Whilst these rules could've been applied earlier, by waiting until after overlap-grouping, unwanted peaks which are harder to filter are often grouped with obviously non-tooth peaks and can therefore be removed safely.

\paragraph{Per Group of Regions Rules:}

A group of overlapping regions is removed if any of the following apply:

\begin{itemize}
    \item Group contains only lingual or buccal pointing surface normals. Any tooth should have both a lingual and a buccal side, or for very slanted teeth, at least a significant variance. The groups on the rugae will all face only palatally so will be rejected by this rule.

    \item Group touches the mesh boundary. This is primarily for intra-oral scans which often pick up bits of cheek which must be ignored.
\end{itemize}

\subsubsection{Group Inline Groups}
\label{sec:group-inline}

This is the second of the two grouping stages.

Labial/buccal and lingual/palatal cusps of molars and premolars will often still be separated but can be put back together by grouping by position along the jaw-line (again allowing indirect grouping). The output groups of groups should be whole teeth.

The arch of the jaw-line makes the geometry of the above deceptively awkward. A quadratic curve, fitted to the horizontal components of the remaining peak points, approximates the jaw-line (see section~\ref{sec:quadratic-curve} for more information on the quadratic). Each region group can be mapped onto the quadratic to determine its span (left-most and right-most position) along the jaw, effectively 1-dimensionalising (flattening) the jaw line. These spans can be compared directly to test if two groups are inline.

\subsubsection{The Output and its Drawbacks}

\begin{figure}[ht]
    \centering
    \captionsetup{margin=.03\textwidth}
    \includegraphics[width=0.9\linewidth]{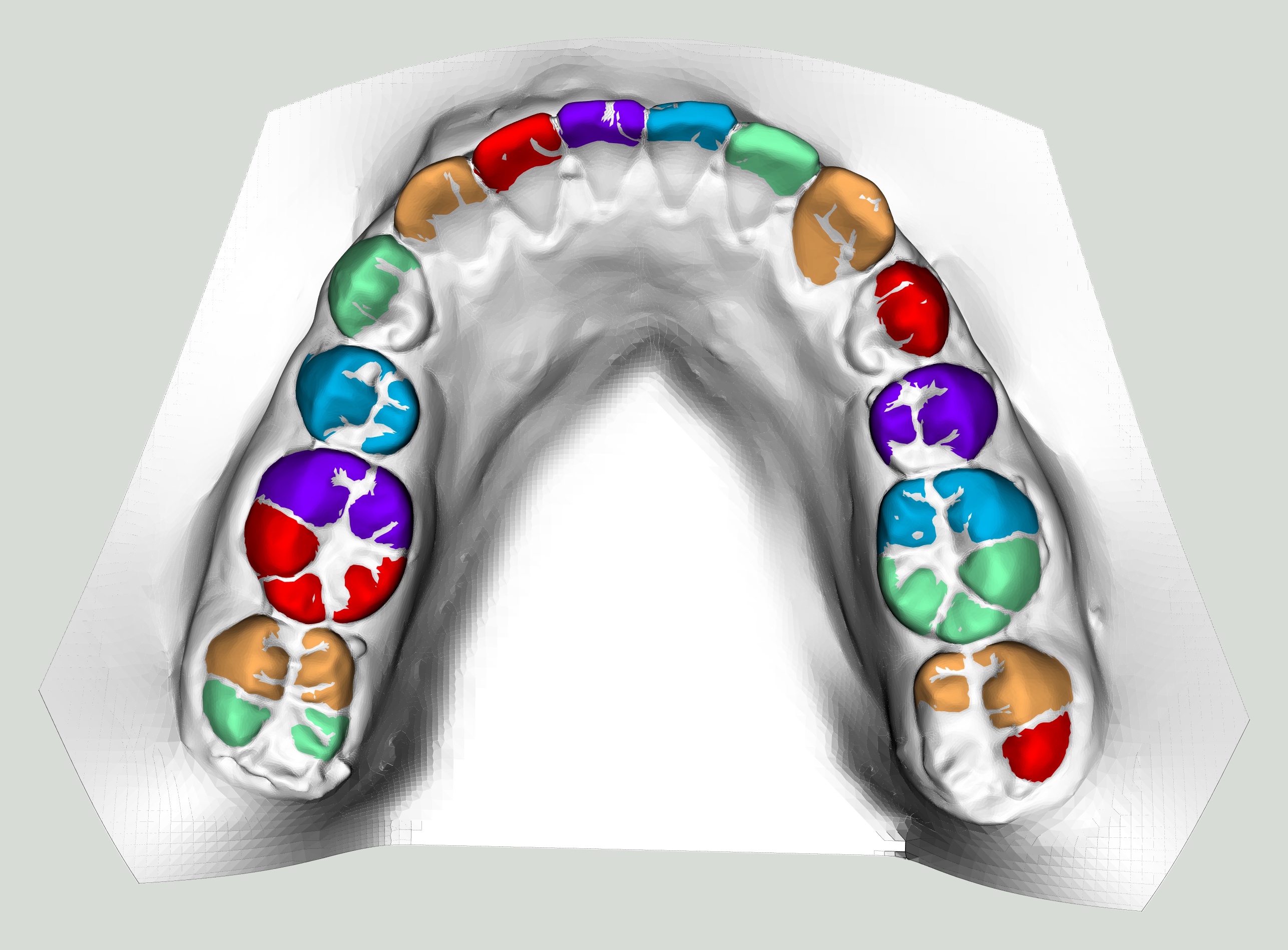}
    \caption{The final output of this tooth segmentation method. Each change in colour signifies a different tooth.}
    \label{fig:partition-output}
\end{figure}

The resulting groups of groups from above are teeth and are shown in figure~\ref{fig:partition-output}.

The only teeth that may still remain split are the mesial and distal halves of the molars (see figure~\ref{fig:partition-output}). These half-molars will need to be merged but we can't do that yet until they are identified as such. Otherwise we may mistakenly glue two premolars together instead. This correction is handled by the tooth assignment step (section~\ref{sec:tooth-assignment}).

Lingual portions of incisors, and the middles and grooves of molars are all inside corners and therefore are unrecognised as parts of the teeth. An adapted version of the method described by Bei-ji Zou e.t al \cite{harmonic-field-ZOU2015132} has been applied as an additional step to capture to whole teeth with good success. But as this has little impact on finding  MHB landmarks, and depends on some heavy-weight sparse linear algebra libraries, applying this step has not been adopted as the default behaviour.

\subsection{Best-fit Quadratic Orientation Curve}
\label{sec:quadratic-curve}

In sections~\ref{sec:partition-into-individual-teeth} and~\ref{sec:tooth-assignment} it is important to be able to perform operations that refer to the arch shape of the model. Namely, to define the directions tangential (mesial and distal) and normal (lingual and buccal) to the jaw's arch, and to facilitate sorting by position around the arch.

To do these requires a continuously differentiable curve fit of the arch. Other works have used cubic-splines to do this \cite{Automatic-Feature-Identification-in-Dental-Meshes,Ma2017-em}. However, it assumes that the points fitted to are ordered and contain no outliers/anomalies (neither assumptions hold here). The jaw is roughly quadratic shaped -- so a simple quadratic fit based on the least squares method was chosen. 

\subsubsection{Construction}
The curve is fitted to the horizontal components of the peak points from section~\ref{sec:peak-points} with $x$-increasing defined as left to right across the mouth and $y$-increasing as going forwards.
\begin{equation}
    y = ax^2 + bx + c
\end{equation}

Most of the time, the fit is tolerant enough that the peak points do not require any filtering beforehand - but not always. It is therefore best to wait until after the peak cleaning from section~\ref{sec:cleaning} has removed those irrelevant points before applying this technique.

\subsubsection{Ordering Peaks and/or Teeth}
\label{sec:sort-peaks}

\begin{figure}[ht]
    \centering
    \captionsetup{margin=.04\textwidth}
	\includegraphics[width=0.9\linewidth]{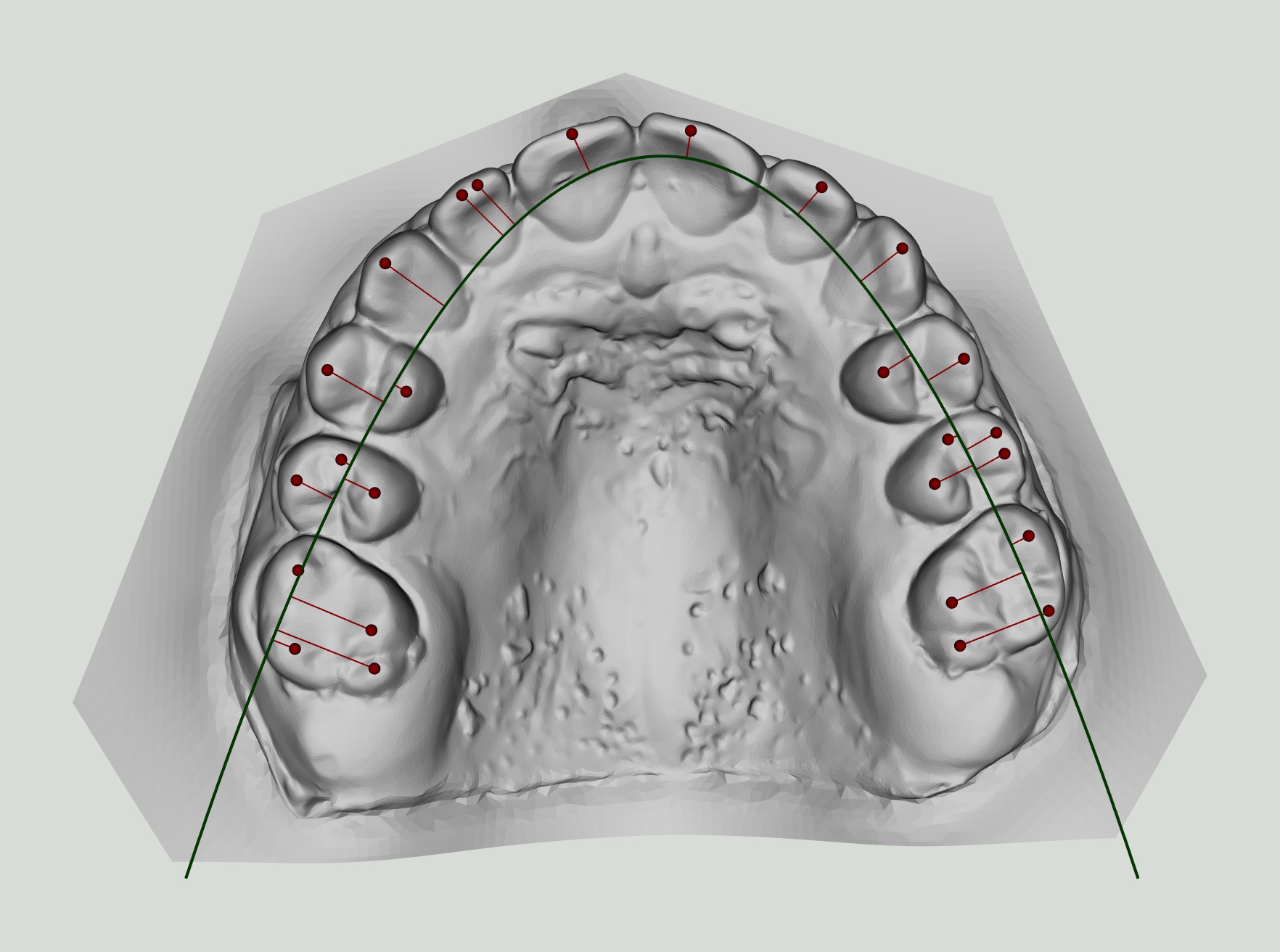}
    \footnotesize{\caption{The least-squares quadratic curve (green line) fitted to the peak points (red dots). Each peak is projected back onto the curve (red lines). The peaks can be ordered by where their projections lie on the curve.}
    \label{fig:qc1}}
\end{figure}

Any series of points, such as the peak points or the centre of mass of each unlabelled tooth, can be objectively sorted and enumerated using the quadratic. To do this, project each point to its nearest point on the curve, then sort and enumerate by the $x$ value of each nearest point (see Figure~\ref{fig:qc1}). 

\subsubsection{Generating \emph{Distal} and \emph{Buccal} Unit-Vectors}
\label{sec:generate-unit-vectors}

\begin{figure*}[!ht]
    \centering
    \begin{subfigure}{0.45\textwidth}
        \captionsetup{margin=.04\textwidth}
    	\includegraphics[width=\textwidth]{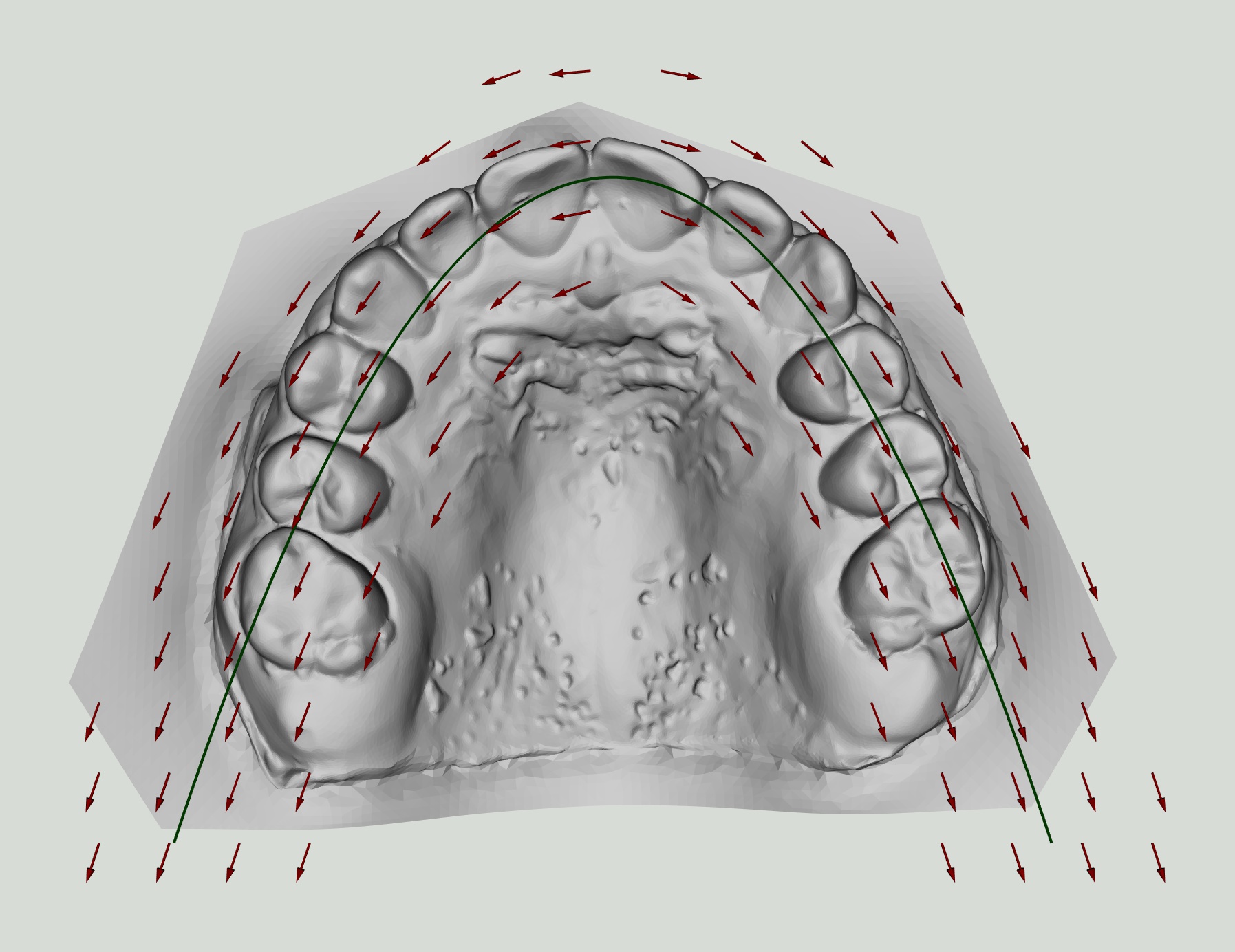}
    	\caption{The arrows indicate the \emph{distal} direction generated by finding tangents to the curve.}
        \label{fig:qc2}
    \end{subfigure} 
    \hspace{.05\textwidth}
    \begin{subfigure}{0.45\textwidth}
        \captionsetup{margin=.04\textwidth}
    	\includegraphics[width=\textwidth]{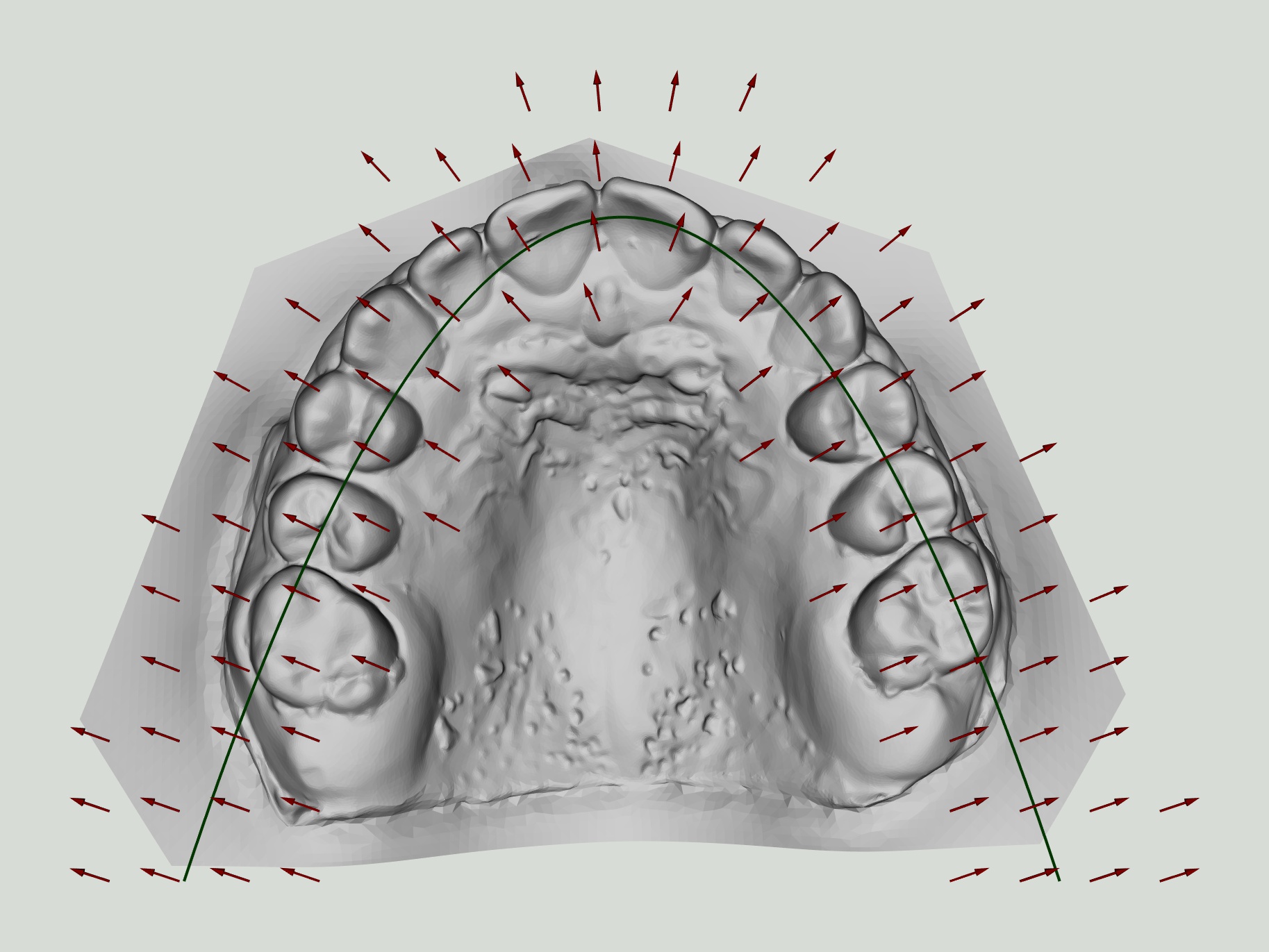}
    	\caption{The arrows indicate the \emph{buccal}/\emph{labial} direction generated by finding outward normals to the curve.}
        \label{fig:qc3}
    \end{subfigure}
    \captionsetup{margin=.025\textwidth}
    \caption{The \emph{distal} and \emph{buccal} directions can be defined at any point in space using the least squares quadratic curve.}
    \label{fig:qc-directions}
\end{figure*}

The directions \textit{along the jaw-line} and \textit{perpendicular to the jaw-line} can be defined with reference to the curve (see figures~\ref{fig:qc2} and~\ref{fig:qc3}). The following commonly required direction vectors can be derived.
\begin{itemize}
    \item \textit{Mesial} (towards the front teeth)  tangent to quadratic with positive y component.
    \item \textit{Distal} (towards the back teeth)  tangent to quadratic with negative y component.
\end{itemize}

\begin{itemize}
    \item \textit{Buccal/labial} (outwards towards the lips)  normal to quadratic with positive y component.
    \item \textit{Lingual/palatal} (inwards towards the tongue)  normal to quadratic with negative y component.
\end{itemize}

The $y$ axis is guaranteed by section~\ref{unit-vector-signs} to point forwards, so a negative $\frac{dy}{dx}$ indicates that the gradient and tangent of the quadratic curve is in the direction of the back of the jaw. As an example, \emph{distal} is calculated in full below.

\begin{equation*}
\text{distal tangent} = 
    \left\lbrace \begin{aligned}
    - \left[ 1, \frac{dy}{dx} \right] \qquad & \text{if } \frac{dy}{dx} \geq 0 \\
    \left[ 1, \frac{dy}{dx} \right] \qquad & \text{if } \frac{dy}{dx} \text{ \textless } 0
    \end{aligned} \right.
\end{equation*}

To convert the direction vector back to 3D use the \emph{right} and \emph{forwards} unit-vectors from section~\ref{sec:orientation} and the $x$ and $y$ components of the vector as follows:

\begin{equation*}
    \text{3D distal} \:\: = \:\: x \: \unitvec{e}_{right} \:\: + \:\: y \: \unitvec{e}_{forwards}
\end{equation*}

At any point in space, the nearest point to the quadratic can be solved for, the tangent or normal be calculated at that nearest point and a unit-vector representing distal (see figure~\ref{fig:qc2}), buccal (see figure~\ref{fig:qc3}), lingual or mesial can be generated.

\subsection{Tooth Assignment}
\label{sec:tooth-assignment}

\begin{figure}[ht]
    \centering
    \captionsetup{margin=.055\linewidth}
    \includegraphics[width=0.9\linewidth]{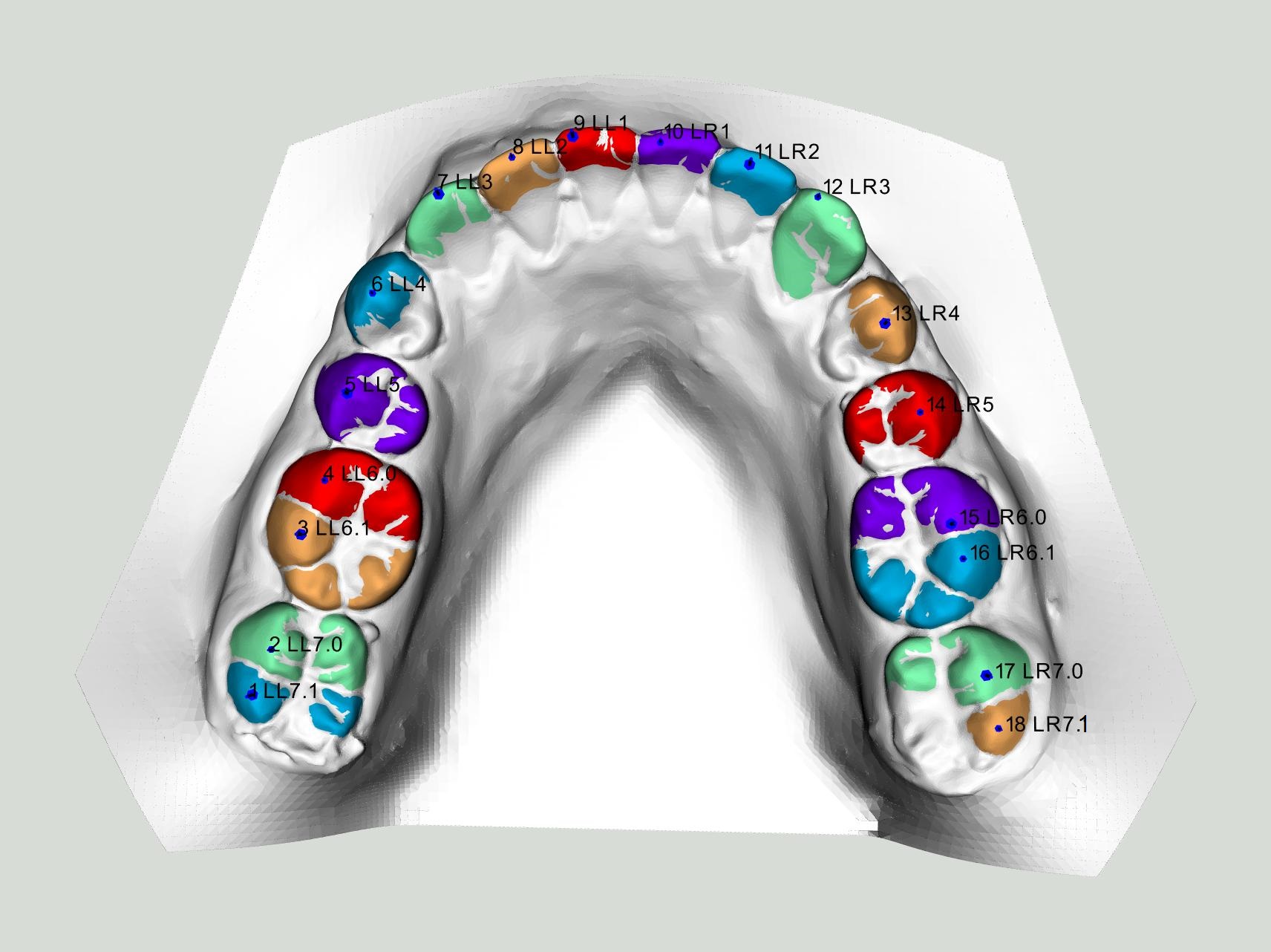}
    \caption{The end goal and final output of the Tooth Assignment method. Each originally unknown tooth from section~\ref{sec:partition-into-individual-teeth} has been labelled with the appropriate tooth type.}
    \label{fig:assignment-goal}
\end{figure}

Section~\ref{sec:partition-into-individual-teeth} yields a collection of unlabelled sub-samples of the original mesh which will be referred to throughout this chapter as \emph{blob}s. Each blob is either a whole tooth, one half of a molar, or occasionally a non-tooth feature. The blobs are sorted and  enumerated (with direction left to right) by section~\ref{sec:sort-peaks}. This section assigns a tooth type to each blob if appropriate. The aimed results are shown by figure~\ref{fig:assignment-goal}. 

Complications that may arise and need dealt with by the method:

\begin{enumerate}
    \item Teeth may be missing, and may or may not leave a space where the tooth would have been.
    \item Molars may be whole, or split into mesial and distal halves, or partially erupted so that only the mesial half is visible.
    \item On adult teeth, the 7s and 8s are much less likely to be present.
    \item There will be occasional non-tooth blobs.
    \item And of course, all the inconsistencies of individual teeth:
    \begin{itemize}
        \item Natural genetic deviation.
        \item Worn, chipped, malformed or re-crowned tips.
        \item Teeth tilted at extreme angles or rotated within their socket.
    \end{itemize}
\end{enumerate}

\subsubsection{Prelude -- Methods That Didn't Work}

Troubleshooting took many unsuccessful attempts until it was possible. A brief history of those attempts is described to demonstrate how the solution was reached.

The first tempting solution is to enumerate round from the centre. This does not work as any missing tooth will lead to wrong count. Furthermore, the centre of the arch isn't precisely known. Absolute position around a jaw line is too inaccurate to be used - an error of a few degrees can lead to picking the wrong tooth. 

Another apparent solution may be to try and use specific tooth characteristics (e.g. number of cusps, geometric properties e.t.c.) to classify teeth. Rule based logic doesn't handle well large varieties of exceptions (or it becomes too complex). The complexity is considerably increased by the range of orientations of teeth. These orientations could be factored out if they were known, but to calculate a tooth's orientation requires that which tooth it is can be known in advance - a \textit{chicken and egg} case.

After failing to create rules for which there would be no exemptions, was realised that to achieve this is unlikely and creates many limitations.
An effective approach should:
\begin{itemize}
    \item Take advantage of general features and characteristics of each tooth type without explicitly relying on them.
    \item Take advantage of neighbouring blobs when considering a particular blob. For example, if blob $B$ looks like either a \nth{1} or \nth{2} right premolar, but the previous blob $A$ is clearly also a premolar then $A$ is the \nth{1} and $B$ is the \nth{2}.
    \item Postpone any \textit{digitisation} (converting from continuous probabilities to a hard yes or no) until as late as possible. 
\end{itemize}

Digitisation should take place only once all observations for each tooth have been combined. Digitisation, in a sense, is an extreme for of rounding and therefore loses information if applied to intermediate results.

With the above conclusions, the following method was constructed.

\subsubsection{Overview}

\begin{figure}[ht]
    \centering
    \captionsetup{margin=.04\linewidth}
    \includegraphics[width=0.9\linewidth]{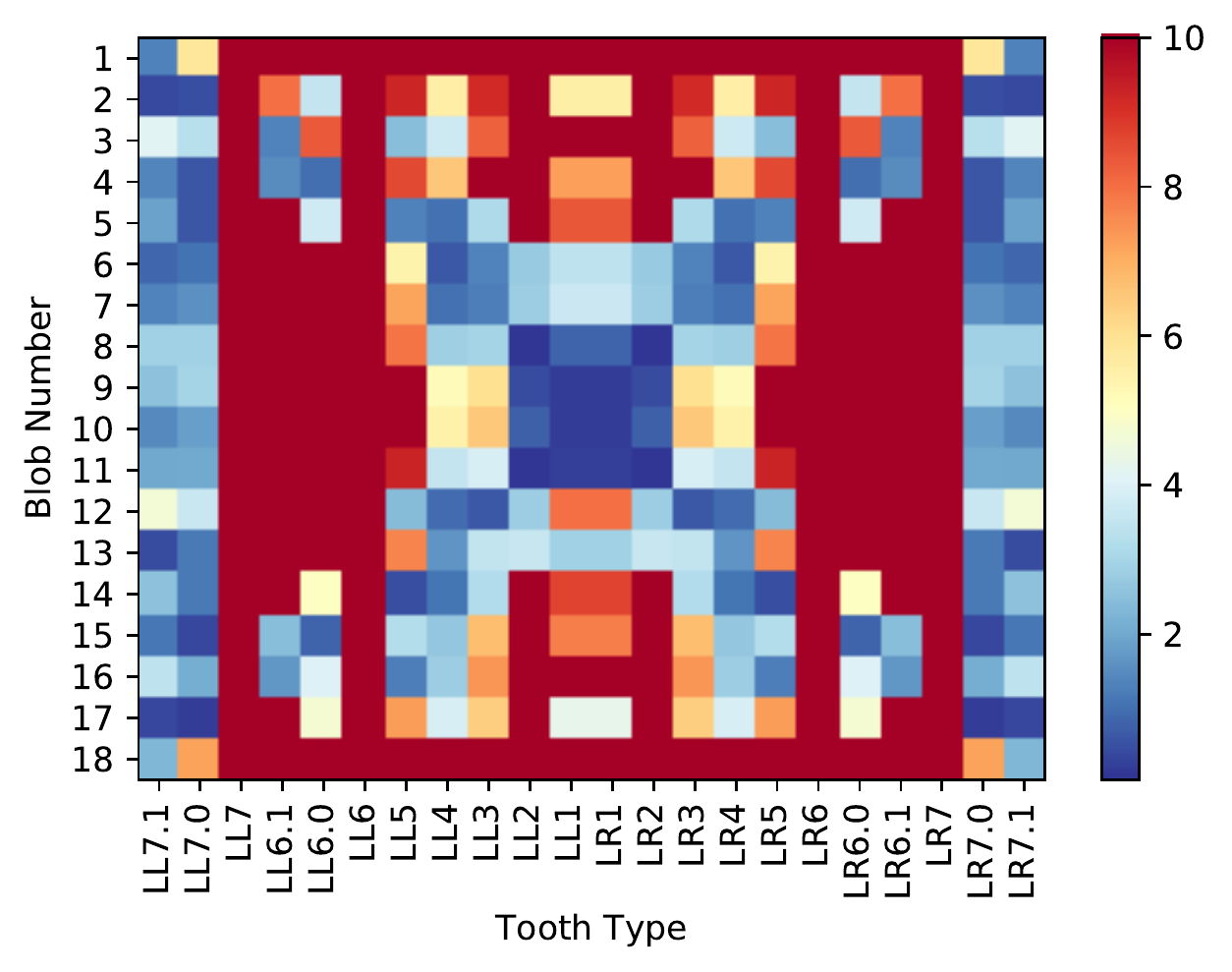}
    \caption{Heat map of assignment costs. Each square represents how unlikely a particular blob (unlabelled tooth-like sub-sample of the mesh) is to be a particular tooth type. Values higher than 10 are clipped to 10 for visual clarity -- costs can be much higher.}
    \label{fig:assignment-costs}
\end{figure}

The method is two-part.

\paragraph{Create a Table of Mismatch Costs:} \emph{Costs} (measures of mismatch) are calculated by comparing each blob to a database of hand-labelled blobs. This is done for all possible blob tooth-type pairs giving a table of costs (such as figure~\ref{fig:assignment-costs}). A mismatch cost is a measure of difference. For example, \emph{the mismatch cost of assigning blob $x$ to LR3 is $0.0$} would mean that blob $x$ is the perfect stereotype of an adult lower right canine. And if the cost was a large number, then blob $x$ would almost certainly not be a LR3. This step will never definitively say that blob $x$ is or isn't a given type -- only how unlikely it is.

\paragraph{Solve the Cost Table:} The table of costs is \textit{solved} to find the \emph{optimal}, defined as lowest possible total mismatch, valid assignment. Thus making the final decision as to which tooth is which. The optimal assignment is found using \emph{Linear Programming}. \footnote{Linear Programming is a wide-spread field of computational mathematics used to minimize (or maximize) a linear objective function subject to linear equality or inequality constraints.}

\paragraph{The tooth types:} searched for are referenced using Palmer notation, e.g. UR2 for upper right \nth{2} incisor, with an extension to describe half-molars: UR6.0 and UR6.1 represent the mesial and distal halves of an UR6. A list of potential tooth types is generated based on the model's dentition type (permanent/deciduous and upper/lower). The tooth types for each jaw type are listed in table~\ref{table:palmers}.

\begin{table}[H]
    \centering
    \tiny
\begin{tabular}{llllllll}
\hline
    \multicolumn{8}{c}{Jaw Types}  \\
\hline
 \multicolumn{4}{c}{Adult} & \multicolumn{4}{c}{Deciduous} \\
 \multicolumn{2}{c}{Upper} & \multicolumn{2}{c}{Lower} & 
 \multicolumn{2}{c}{Upper} & \multicolumn{2}{c}{Lower} \\
\hline
 UL1   & UR1   & LL1   & LR1   & ULA   & URA   & LLA   & LRA   \\
 UL2   & UR2   & LL2   & LR2   & ULB   & URB   & LLB   & LRB   \\
 UL3   & UR3   & LL3   & LR3   & ULC   & URC   & LLC   & LRC   \\
 UL4   & UR4   & LL4   & LR4   & ULD   & URD   & LLD   & LRD   \\
 UL5   & UR5   & LL5   & LR5   & ULE   & URE   & LLE   & LRE   \\
 UL6   & UR6   & LL6   & LR6   & ULE.0 & URE.0 & LLE.0 & LRE.0 \\
 UL6.0 & UR6.0 & LL6.0 & LR6.0 & ULE.1 & URE.1 & LLE.1 & LRE.1 \\
 UL6.1 & UR6.1 & LL6.1 & LR6.1 &       &       &       &       \\
 UL7   & UR7   & LL7   & LR7   &       &       &       &       \\
 UL7.0 & UR7.0 & LL7.0 & LR7.0 &       &       &       &       \\
 UL7.1 & UR7.1 & LL7.1 & LR7.1 &       &       &       &       \\
 UL8   & UR8   & LL8   & LR8   &       &       &       &       \\
 UL8.0 & UR8.0 & LL8.0 & LR8.0 &       &       &       &       \\
 UL8.1 & UR8.1 & LL8.1 & LR8.1 &       &       &       &       \\
\hline
\end{tabular}
\caption{All tooth types searched for, for each jaw type.}
\label{table:palmers}
\end{table}

Some clarifications to make here:

\begin{itemize}
    \item Equivalent teeth from different jaw types such as \{UR1, LR1, URA, LRA\} are treated independently as if they weren't related. 
    \item Tooth sub-types (incisors, canines, premolars, molars) are similarly ignored. An UR1 is (to the software's mind) unrelated to an UR2.
    \item Both halves of a molar as well as the whole molar are thought of as independent whole teeth throughout most of this method. The \textit{solve the cost table} stage adds additional constraints to the linear programming model stating that a whole molar and its two halves are mutually exclusive.
    \item Whilst creating the cost table, left tooth types are considered equivalent to their corresponding right types (making the cost table symmetrical), but whilst solving it, they are searched for separately.
\end{itemize}

\subsubsection{Create a Cost Table}
\label{sec:create-cost-table}

An unknown-blob tooth-type pair is tested by comparing the blob to a database of hand labelled blobs (training set). Rather than try to compare whole blob meshes directly (which was tried unsuccessfully), measurements of the meshes are taken and compared. These measurements are referred to in this paper as \emph{tooth characteristics}. They must be applicable to any tooth type and yield a single number per blob.

\paragraph{Tooth Characteristics}
\label{sec:tooth-characteristics}

\begin{figure}[ht]
    \centering      
    \captionsetup{margin=.04\textwidth}
    \includegraphics[width=0.9\linewidth]{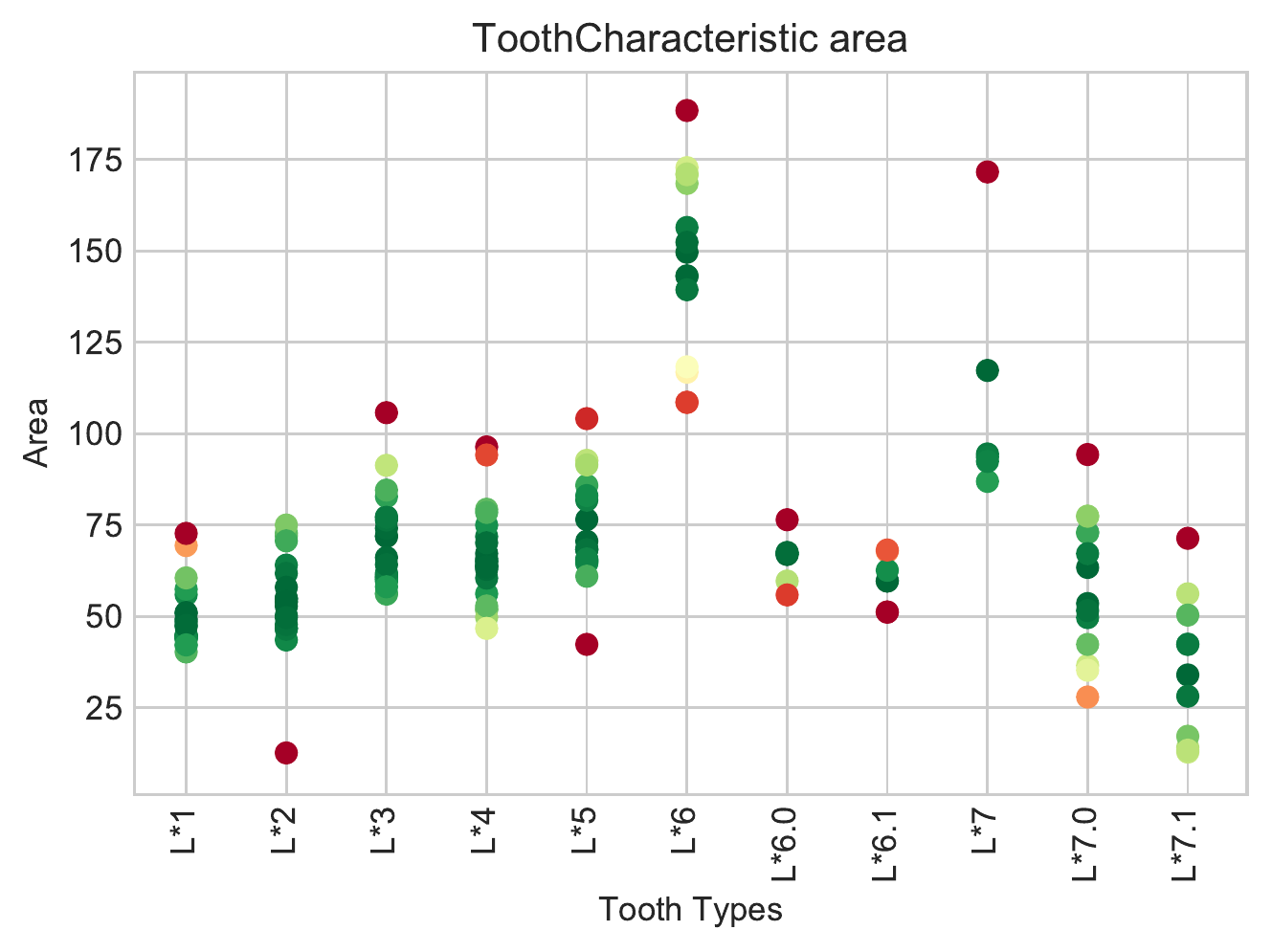}
    \caption{The surface area of each tooth in the training set, grouped by tooth-type. The (*) placeholder in the tooth-type names symbolises that left and right teeth can be used together.}
    \label{fig:tooth-characteristic}
\end{figure}

 The simplest example is the Area tooth characteristic which measures the total surface area in mm$^2$ of a blob. The total surface area has been measured on all the hand-labelled blobs in the training set and the results were grouped by tooth type (see figure~\ref{fig:tooth-characteristic}). These values are referred to as \emph{reference values}. When analysing an unknown blob, the surface area is calculated for that blob (the \emph{test value}) and compared to each group of reference values. The more the test value differs from the reference values, the higher the cost for that tooth type.

To demonstrate, suppose an unknown blob had a surface area of $75\text{mm}^2$. 

\begin{itemize}
    \item $75$ is an ideal value for the blob to be a 3, 4 or 5, or the mesial half of a 7. Costs should be small, typically $<1.0$.
    \item $75$ is on the high side for the blob to be either of the incisors or the distal halves of the molars. Costs should be medium, typically $2.5 - 5.0$.
    \item $75$ is far too low to be either whole molar. Costs should be very high, typically $10.0 - 100.0$.
\end{itemize}

This is done for each unknown blob, for each tooth characteristic. The costs for each characteristic (listed in table~\ref{sec:tooth-characteristics}) are averaged to give an overall cost for each blob tooth-type pair.

This design makes it possible to provide overall, qualitative trends about tooth types without turning them into strict rules which must apply to every tooth. For example, one might say \emph{canines and premolars have pointed tips but incisors have flat tops}. On its own this statement is unhelpful to a computer because it's not always true. A worn canine is not pointy and would confuse a \emph{pontiness} rule into thinking it's an incisor. But turning it into tooth characteristic overcomes the potential complication. First, a method to quantify a blob's pointiness must be devised. The tooth characteristic will then handle everything else. It will:

\begin{itemize}
  \item Train itself by applying the method to all labelled blobs in the training set to understand exactly how pointy each tooth type should be. 
  \item Handle the comparison of test values with reference values.
  \item Self evaluate its reliability. A tooth may be worn, so pointiness is not wholly indicative. Through use of the error metric (section~\ref{sec:error-metric}), less trustworthy characteristics are given less weight in the final outcome.
\end{itemize}

There is no limit to how many characteristics can be used. Generally, the more the better, but it transpired that not many are necessary. The list chosen is given in table~\ref{table:tooth-characteristics}. 

\begin{table}[H]
    \centering
    \scriptsize
\begin{tabular}{|p{.16\linewidth}|p{.34\linewidth}|p{.34\linewidth}|}
\hline
 Name & Description & Goal  \\
\hline
 Total Surface Area           & Calculates the total surface area of the blob. This is very easy to implement reliably.      &           \\
 Mesiodistal Width            & As the name suggests. This requires a \emph{distal} unit-vector which can be generated as described in section~\ref{sec:generate-unit-vectors}.      &         \\
 Buccolingual Width           & Similarly, this requires a \emph{buccal} unit-vector &  This separates premolars and molars from incisors and canines far more reliably than counting cusps. \\
 Pointiness                   & Calculated by measuring the mesiodistal width at 1mm from the top and comparing it to the mesiodistal width from above.  &  Helps to distinguish canines from incisors.     \\
\hline
\end{tabular}    
\caption{The tooth characteristics currently used.}
\label{table:tooth-characteristics}
\end{table}

A few other characteristics that didn't work are noted in table~\ref{table:reject-tooth-characteristics} to deter anyone from trying them again. The results were too noisy to be of use and in some cases were very processor demanding.

\begin{table}[H]
    \centering
    \scriptsize
\begin{tabular}{|p{.15\linewidth}|p{.35\linewidth}|p{.34\linewidth}|}
\hline
 Name & Description & Goal  \\
\hline
 Symmetry           & Mirrors the blob them aligns and compares the mirrored version to the original using Iterative Closet Point Algorithm.   &   Premolars are very symmetric mesiodistally whereas half-molars aren't.   \\
 Total Curvature    & Add all the unsigned curvatures together.  &  Supposed to rather lazily separate based on how textured each tooth is.       \\
\hline
\end{tabular}    
\caption{Some old, less successful tooth characteristics.}
\label{table:reject-tooth-characteristics}
\end{table}

\paragraph{Cost Metric -- Comparison Function}
\label{sec:error-metric}

The cost metric quantifies the mismatch between a test value and a set of reference values (for a single characteristic). The main reason for this choice of metric was to avoid any arbitrary weights that have to be machine-learnt. 

The metric is based on the mean square error (MSE), which is a common default in machine learning. For a test value $t$ and vector $\mathbf{r}$ of length $k$ of reference values:

\begin{equation*}
    MSE\:(t,\: \mathbf{r}) \quad := \quad \frac{1}{k} \: \sum_{i=0}^k \: (t - r_i)^2
\end{equation*}

In order to be able to meaningfully combine and compare costs and to better solve this optimal assignment problem using Linear Programming, a couple of modifications are required to conform to the following requirements. For a given set of reference values:

\begin{enumerate}
    \item The minimum possible cost should be zero. With MSE, it depends on the reference values.
    \item The output costs should be scale independent. With MSE, a characteristic with large typical values will dominate another with smaller ones.
    \item High inter tooth-type variance in the reference values should reduce the costs for test values, thus making the characteristic more lenient.
\end{enumerate}

Requirement 1 is trivial to achieve. Just find the minimum possible cost $C_{min}$ and subtract it from future calculated costs. Both 2 and 3 can be achieved simultaneously by feeding the metric each of the reference values as test values, taking mean $\overline{C_{ref}}$ of the resulting costs and dividing through any future costs by this mean. The final formula is written below:

\begin{equation*}
    C(t) \: = \: \frac{MSE\:(t,\: \mathbf{r}) \: - \: C_{min}} {\overline{C_{ref}}}
\end{equation*}


\subsubsection{Solve the Cost Table} 
\label{sec:solve-cost-table}

\begin{figure}
    \centering
    \captionsetup{margin=.06\linewidth}
    \includegraphics[width=0.9\linewidth]{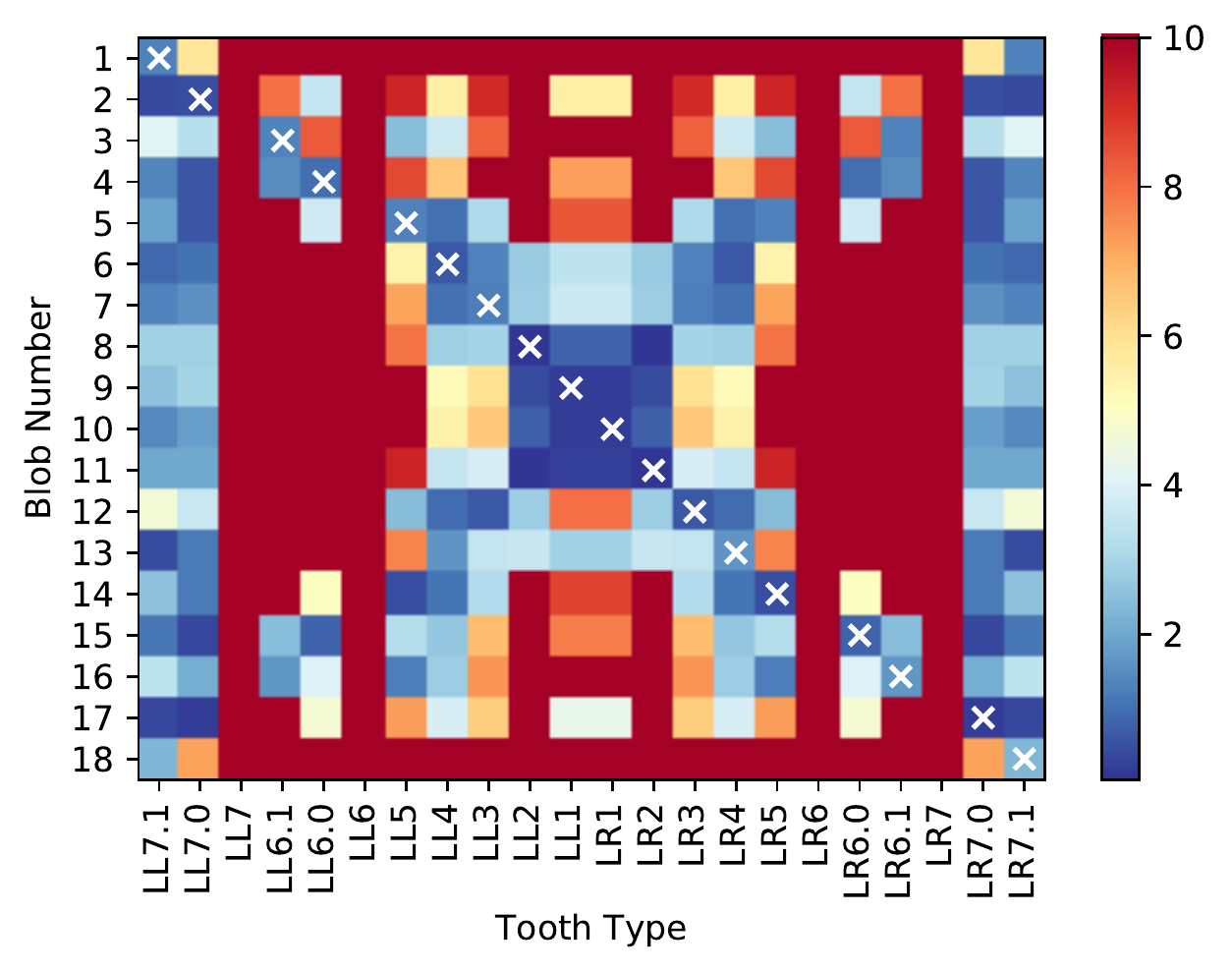}
    \caption{The heatmap from figure~\ref{fig:assignment-costs} with the chosen \emph{optimal} assignment (white X markers). The X markers always form a downward diagonal to ensure the tooth-types are in order along the jaw and will aim to be only on dark blue (low costing) squares. This is from the same model as figure~\ref{fig:assignment-goal} so the X markers here match the annotations there.}
    \label{fig:assignments}
\end{figure}

The lowest costing possible assignment that is \emph{valid} is found using Linear Programming. General purpose linear programming solver packages are freely available. This implementation uses \href{https://coin-or.github.io/pulp/index.html}{PuLP}. All that is required is to formulate the problem in the mathematical format that is standard in linear programmings. The results are shown in figure~\ref{fig:assignments}. LP problems can easier to build from scratch than to read but the gory details of the LP problem are included below.

\paragraph{Decision Variables}

Assignments are represented with Boolean (true or false) variables. Let $D$ be an $m \times n$ table of Booleans, where $m$ is the number of blobs and $n$ is the number of tooth types.

\begin{equation*}
    D[i, j] = 
    \left\lbrace \begin{aligned}
        1,& \quad\text{If the $i^{th}$ blob is of the $j^{th}$ tooth type.} \\
        0,& \quad\text{Otherwise.}
    \end{aligned} \right.
\end{equation*}

A blob may not be a tooth. And a tooth type may be missing. These need to be accounted for. 

\begin{align*}
    NT[i] &= \left\lbrace \begin{aligned}
        1,& \quad\text{If the $i^{th}$ blob is non-tooth.} \\
        0,& \quad\text{Otherwise.} 
    \end{aligned} \right. \\
    MI[j] &= \left\lbrace \begin{aligned}
        1,& \quad\text{If the $j^{th}$ tooth-type is missing.} \\
        0,& \quad\text{Otherwise.} 
    \end{aligned} \right. \\
\end{align*}

\paragraph{Constraints} 

All the rules of dentistry must be expressed as mathematical constraints. 

Each blob can either be exactly one tooth type, or it could be non-tooth. So for each $i = \lbrace 1 \dots m \rbrace$:

\begin{equation*}
    \sum_{j=1}^n D[i, j] + NT[i] = 1
\end{equation*}

Each tooth type can either appear exactly once or that tooth type is missing: For each $j = \lbrace 1 \dots n \rbrace$:

\begin{equation*}
    \sum_{i=1}^n D[i, j] + MI[j] = 1
\end{equation*}

Tooth types should appear in the correct order. To express this mathematically we utilize the fact that both the tooth types and the blobs are ordered from left to right along the jaw. To enforce the order we then need only to ensure that the assignment table $D$ does not reorder them. If blob $i$ is of tooth type $j$, then no previous blobs can be of further right tooth types and no later blobs of further left tooth types. So if $D[i, j]$ then none of,
\begin{align*}
    D[i', j'], \quad  i'= \lbrace 1, \dots, i-1 \rbrace, \:\: j'= \lbrace j+1, \dots, n \rbrace \\
    D[i', j'], \quad  i'= \lbrace i+1, \dots, m \rbrace, \:\: j'= \lbrace 1, \dots, j-1 \rbrace
\end{align*}

This mutual exclusiveness can be expressed mathematically by adding the variables and applying an upper bound ($\leq$) to the sum.

\begin{align*}
        D[i, j] \:\: + \:\: \frac{1}{m \times n} \:\: \sum_{i' = 1}^{i-1} \:\: \sum_{j' = j+1}^{n} \:\: D[i', j'] \quad & \leq \quad 1 \\
        D[i, j] \:\: + \:\: \frac{1}{m \times n} \:\: \sum_{i' = i+1}^m \:\: \sum_{j' = 1}^{j-1} D[i', j']  \quad & \leq \quad 1
\end{align*}

The $\frac{1}{m \times n}$ is necessary so that multiple elements inside the double $\Sigma$ sum can be true simultaneously without violating this constraint.

A whole molar must be missing if either of it's half molar types are not missing: For each molar,

\begin{equation*}
    MI[\: \text{whole} \:] + \frac{1}{2}(\: MI[\: \text{mesial} \:] + MI[\: \text{distal} \:]\:) \:\: \geq \:\: 1    
    \label{eq:half-molar-constraint}
\end{equation*}


Note the use of $\geq$ instead of $=$. This is because if a patient was missing the molar then the LHS of the equation above would be $2$. Similarly, if the molar were partially erupted then the LHS would be $1.5$.

\paragraph{Objective Function} 

The objective is to minimise the total mismatch cost which is the sum of each entry in the cost table whose corresponding entry in the assignment table is true (left half of equation~\ref{eq:objective-function} below). 

There also has to be some form of penalty for marking blobs as non-tooth and/or marking tooth-types as missing. Without such a penalty, the optimal solution, with a net cost of $0$, would unconditionally be to mark all the blobs as no-tooth and all teeth as missing. Penalising for missing tooth-types gives the advantage of being able to vary the penalty for different tooth-types (right half of equation~\ref{eq:objective-function}). This provides a very convenient way to tell Linear Programming that 7s and 8s (wisdom teeth) are significantly less common than the other tooth types. $P[j]$ from above represents the probability that a patient will have a $j^{th}$ tooth-type, which can be trivially derived from the training data.

\noindent Minimise:
\begin{equation}
    \quad \sum_{i=i}^m \:\: \sum_{j=1}^n \:\: D[i, j] \times C[i, j] \:\:\: 
     + \:\:\: 8.0 \: \sum_{j=1}^n \:\: P[j] \times MI[j]
     \label{eq:objective-function}
\end{equation}

The penalty multiplier (or \textit{fussiness factor}) $8.0$ is arbitrary. It controls how atypical a tooth can be before it is assumed to be non-tooth.

\section{Results}

\newlength{\figwidth}
\setlength{\figwidth}{.8\linewidth}


This section serves as a graphical results section and highlights some of the more prominent issues either tackled or still to resolve. Throughout this section, all models are coloured by the output from the tooth assignment step (i.e. one colour per tooth, no colour for an unrecognised tooth) and marked with an annotated black cross-hairs on each MHB landmark.

\subsection{Orientation}

\begin{figure}[ht]
    \centering
    \captionsetup{width=\figwidth}
    \includegraphics[width=\figwidth]{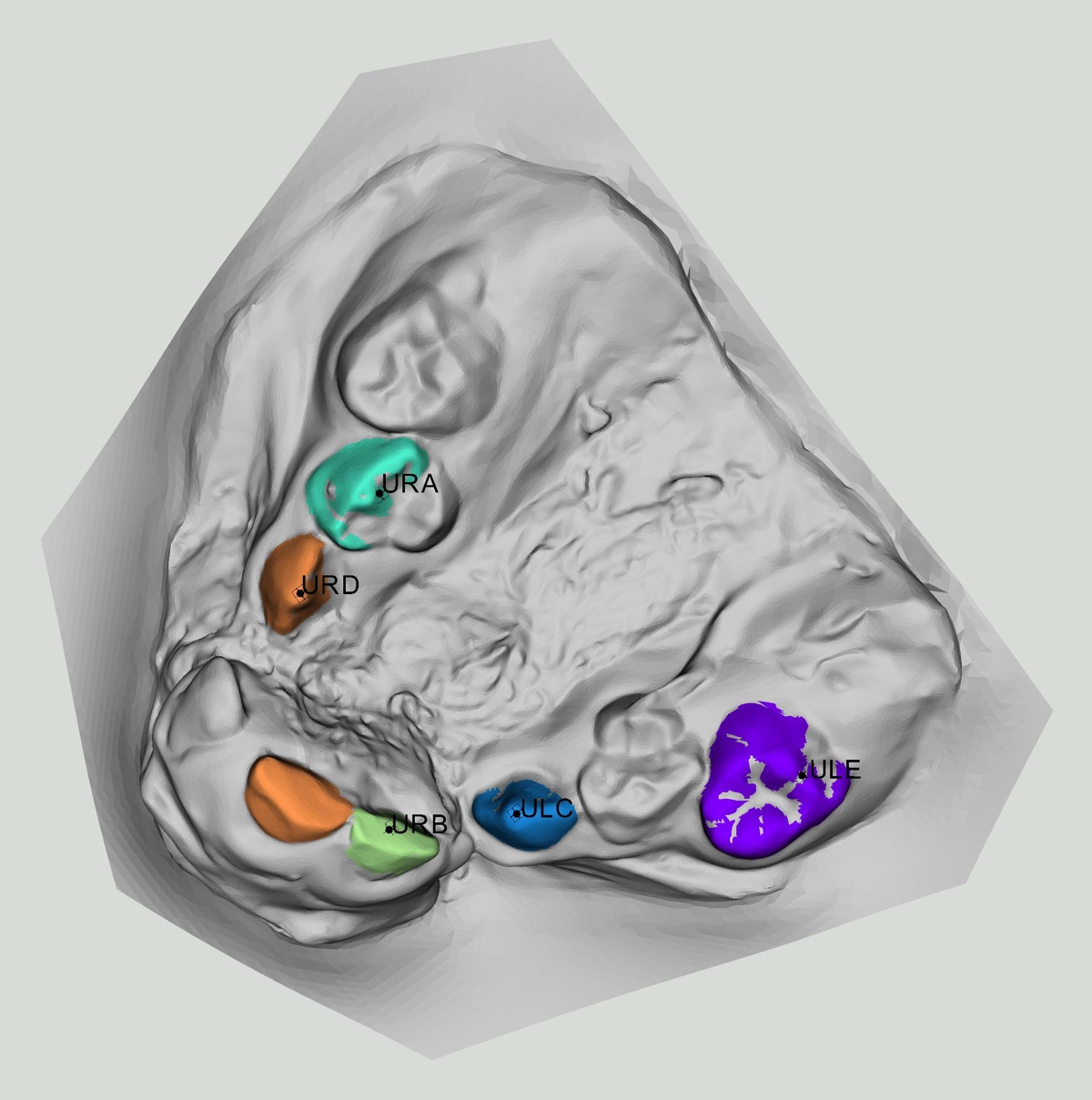}
    \caption{The only model that failed at the orientation step. This photo is oriented with what the software mistakenly thought was the front of the model at the top of the image. Of course, with the orientation wrong, every subsequent output is nonsense.}
    \label{fig:pca-gone-bad}
\end{figure}

Orientation is easily the most reliable step despite its having the least information to work with. There was only one model, shown in figure~\ref{fig:pca-gone-bad}, for which it didn't work. In this case the cause was a double-cleft which is so heavily textured (giving it a higher vertex density than the teeth) that it dominates PCA. 

\subsection{Tooth Partitioning}

\begin{figure}[ht]
    \centering
    \captionsetup{width=\figwidth}
    \includegraphics[width=\figwidth]{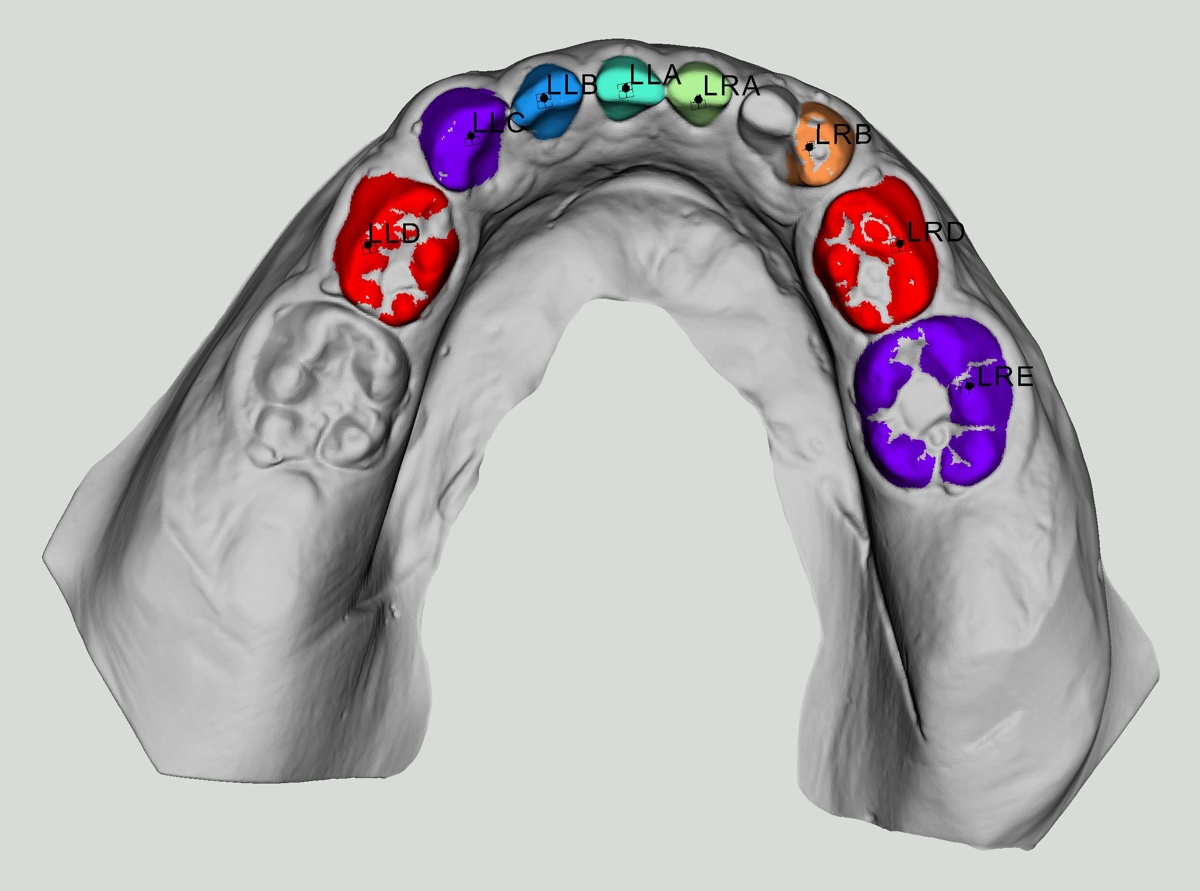}
    \caption{Model with both flattened incisors and a weak outline around the LLE resulting in the LLE not being recognised.}
    \label{fig:partition-error-0}
\end{figure}

Tooth partitioning (section~\ref{sec:partition-into-individual-teeth}) gave mixed degrees of tolerance to unclear boundaries and rough tooth surfaces. The adaptive curvature threshold allows it to use a low threshold for models with weak outlines or a high threshold for models with noisy/bumpy tooth surfaces. However, it can not do both simultaneously so a model with at least one poor outline and one bumpy tooth surface will always lose at least one of them. Flattened incisor tips (very common) form dimples at the top giving the same affect as bumpy teeth i.e. it forces the the curvature threshold up so that the algorithm can cross the dimples to reach the rest of the tooth. Figure~\ref{fig:partition-error-0} demonstrates exactly this case. This is the most common error that the software makes.

Perhaps a more resilient algorithm in the future will be able to choose thresholds per-tooth rather than per-model. Our attempts to do so generally compromised the cleaning steps in section~\ref{sec:cleaning}. 

\begin{figure}[ht]
    \centering
    \captionsetup{width=\figwidth}
    \includegraphics[width=\figwidth]{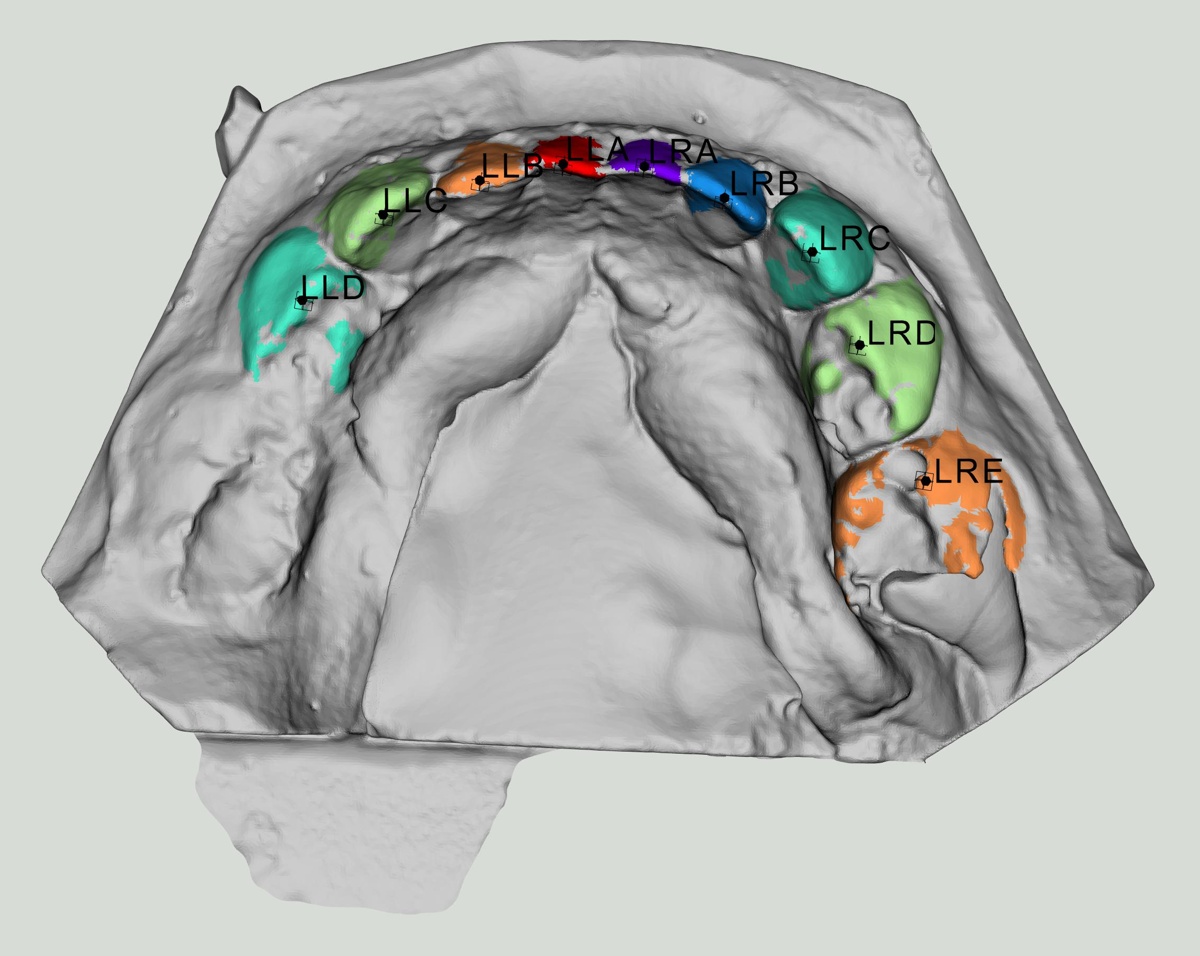}
    \caption{A very poor impression of an infant.}
    \label{fig:partition-error-1}
\end{figure}

The software exceeded expectations on some really poor quality casts (see figure~\ref{fig:partition-error-1}).

\begin{figure}[ht]
    \centering
    \captionsetup{width=\figwidth}
    \includegraphics[width=\figwidth]{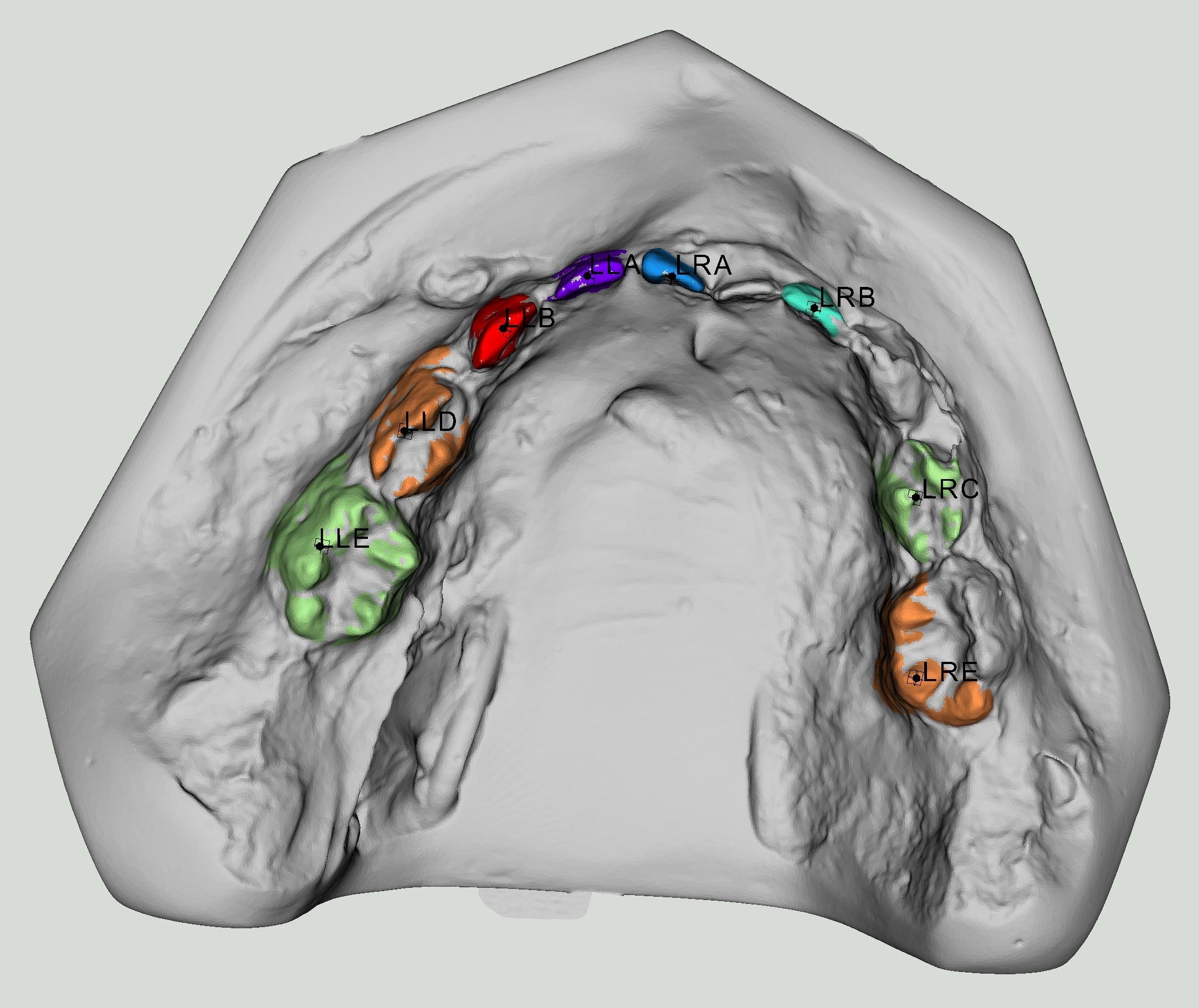}
    \caption{Our very worst plaster model. Whilst the software has found more teeth than we might expect, it is very difficult to find landmarks on teeth that are as poorly recognised as the LRE.}
    \label{fig:partition-error-2}
\end{figure}

Its tolerance is of course finite. See figure~\ref{fig:partition-error-2}). As a side-note: Perhaps, before any further work is done on automating their interpretation, some investigation into getting better impressions/scans from infants would be appropriate.

\subsubsection{Crowding}

\begin{figure}[ht]
    \centering
    \captionsetup{width=\figwidth}
    \includegraphics[width=\figwidth]{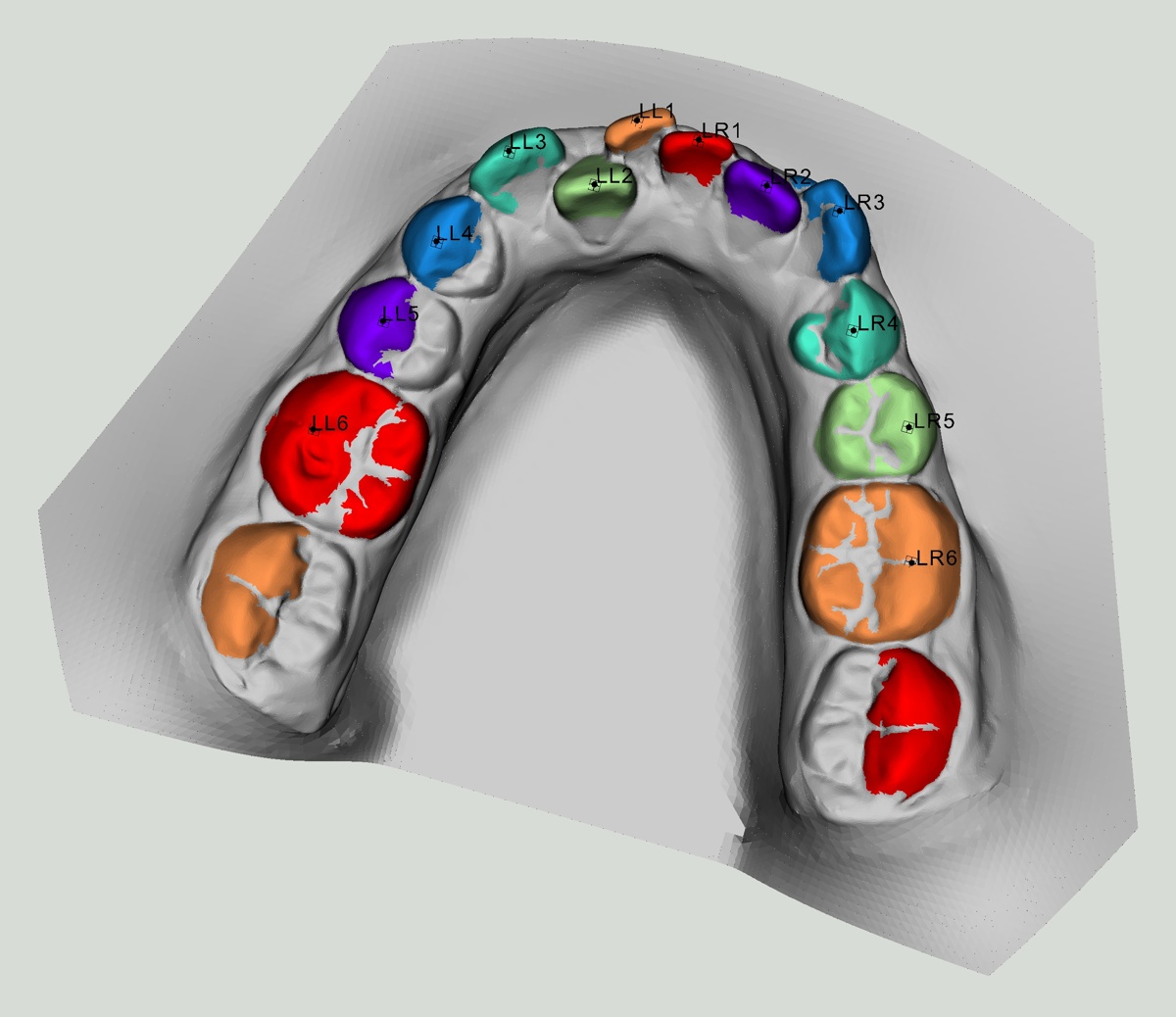}
    \caption{A model with some crowding. The software is tolerant to this degree of crowding (but not much more).}
    \label{fig:crowding-ok}
\end{figure}

\begin{figure}[ht]
    \centering
    \captionsetup{width=\figwidth}
    \includegraphics[width=\figwidth]{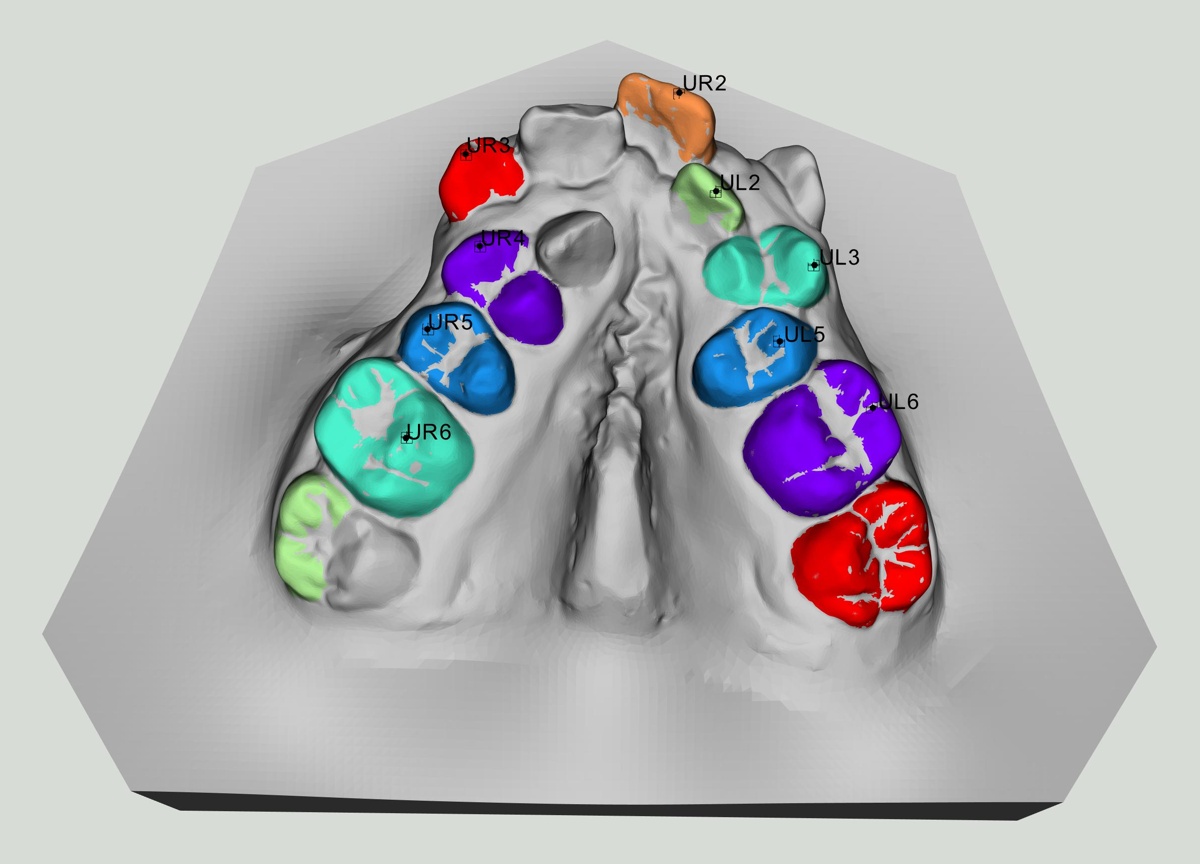}
    \caption{A model with too much crowding for the software. The position of UL2 behind the UL1 directly contradicts the assumption, made in section~\ref{sec:group-inline}, that such a positioning implies that those two teeth must in-fact be two parts of the same tooth e.g two cusps of a pre-molar. After tooth partitioning mistakenly decides that the UR1 and UR2 are one tooth, tooth assignment then, on failing to find a tooth type that resemble this strange \emph{double incisor} tooth, rejects it as non-tooth.}
    \label{fig:too-much-crowding}
\end{figure}

It took some coaxing of overlap thresholds in section~\ref{sec:group-inline} but we were able to get good performance for tooth crowding. The difficulty is that lingual and buccal halves of premolars and molars are paired up only because they are in the same position around the jawline so two incisors that are sufficiently crowded together will be mistakenly grouped. We've managed to give it enough tolerance to allow cases like the one shown in figure~\ref{fig:crowding-ok} but, by design, this algorithm will always fail for cases such the one in figure~\ref{fig:too-much-crowding}.

\subsubsection{Intra-oral scans}

Intra-oral scans impose quite a different set of problems. The tooth outlines are, provided they are scanned properly, much crisper but the tooth surfaces often obtain a \emph{fuzzy} texture which negates the advantage of clear outlines. It is primarily for these models that the adaptive curvature threshold was needed. The chosen threshold is typically much much higher for intra-oral scans.

\begin{figure}[ht]
    \centering
    \captionsetup{width=\figwidth}
    \includegraphics[width=\figwidth]{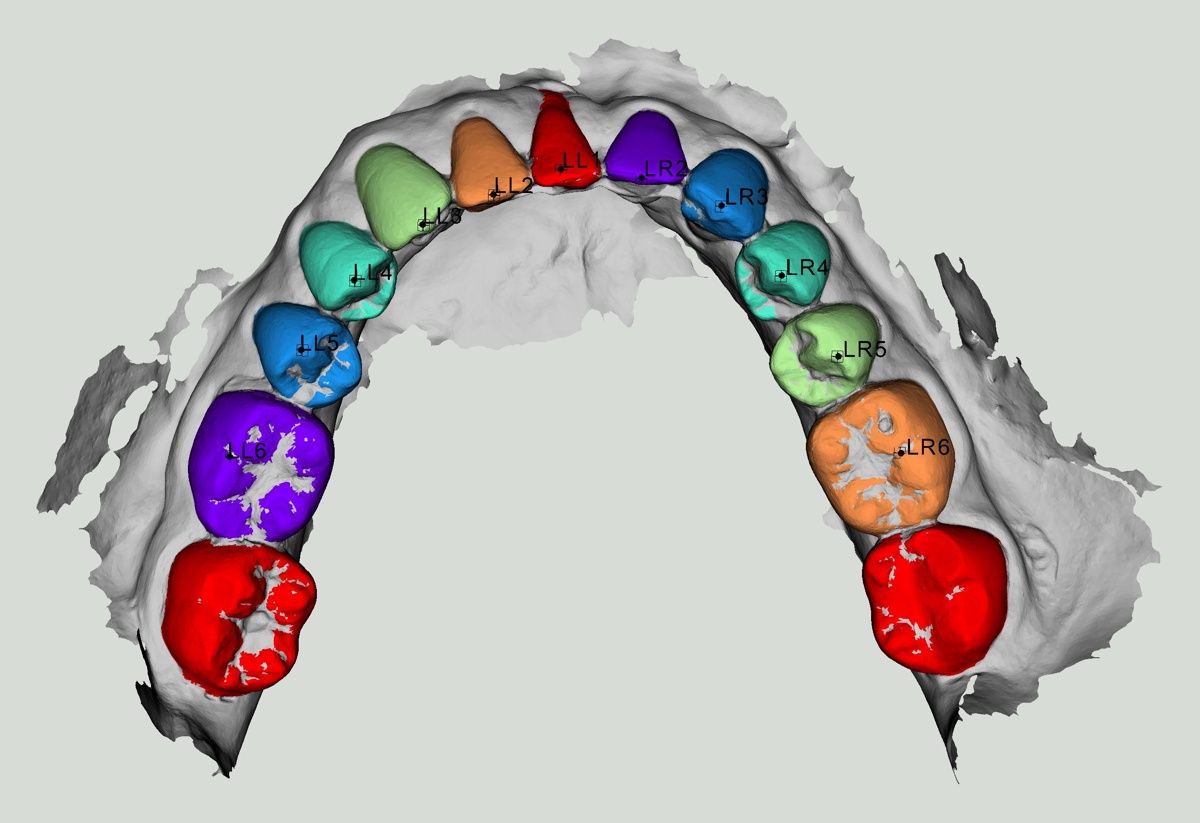}
    \caption{An intra-oral scan with somewhat chaotic trimming. The software is able to correctly ignore it.}
    \label{fig:happy-io}
\end{figure}

\begin{figure}[ht]
    \centering
    \captionsetup{width=\figwidth}
    \includegraphics[width=\figwidth]{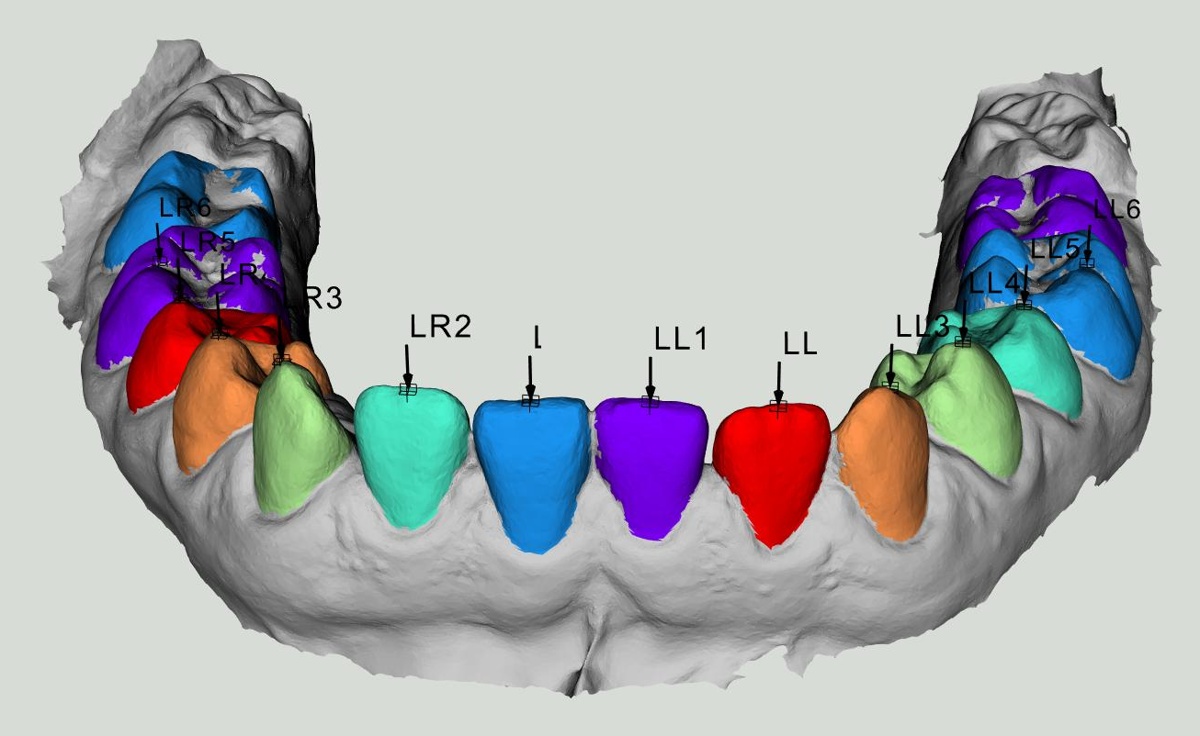}
    \caption{An nice intra-oral model. The scan includes a non-zero amount of gum surrounding each tooth. Partitioning works OK.}
    \label{fig:uncropped}
\end{figure}

\begin{figure}[ht]
    \centering
    \captionsetup{width=\figwidth}
    \includegraphics[width=\figwidth]{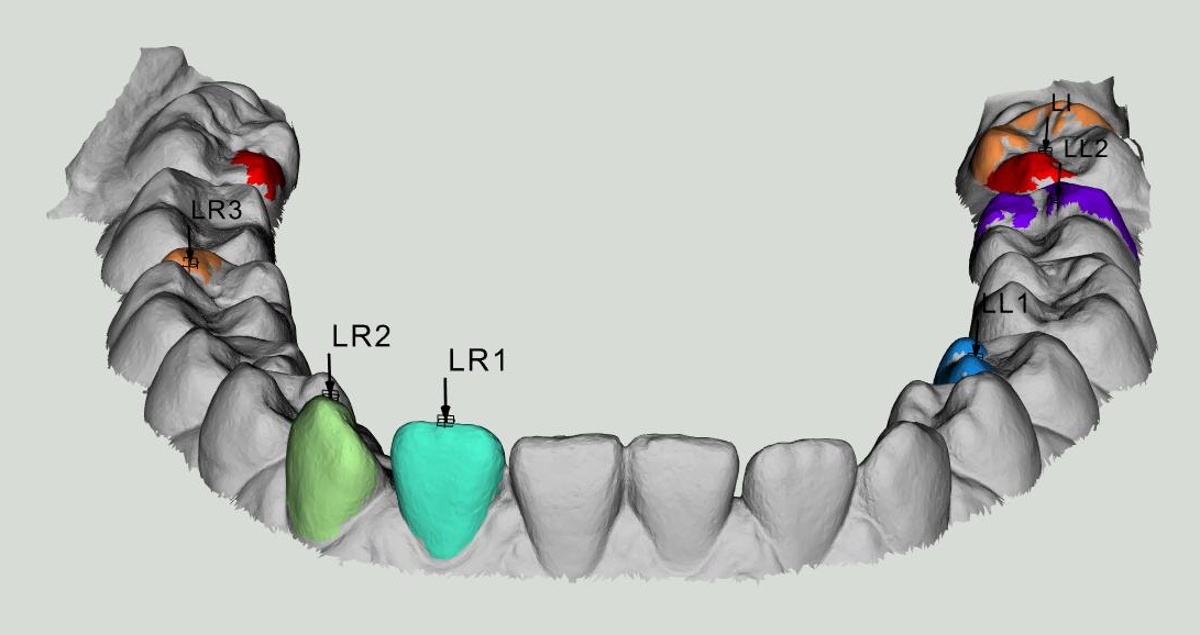}
    \caption{The same model as shown in figure~\ref{fig:uncropped} but with the bottom cropped off. Any tooth which now touched the edge of the scan is rejected in section~\ref{sec:cleaning}.}
    \label{fig:cropped}
\end{figure}

Another problem is the \emph{trimming} or where the scan stops. Partially scanned bits of cheek, lips or tongue often appear in these scans and collect peak points which need to be ignored. A simple \emph{reject anything that touches the mesh boundary} rule (section~\ref{sec:cleaning}) easily gets rid of them as in figure~\ref{fig:happy-io}. However, there is a downside to this -- by definition any \textit{tooth} which touches the mesh boundary is lost. Thus, with this rule in effect, clinicians are required to scan to the base any teeth they wish to analyse. To see the effect of this compare figure~\ref{fig:cropped} to figure~\ref{fig:uncropped}.

In our dataset of intra-oral scans, it was rare that any incisors, canines or premolars weren't fully scanned but about 15\% of molars were lost because of this. To fully capture a molar requires getting the scanner behind it to capture its distal side -- an uncomfortable procedure. This also has the potential to waste a lot of clinician's time should a patient need to be rescanned because of a small gap in a tooth capture.

\subsection{Tooth assignment}

\begin{figure}[ht]
    \centering
    \captionsetup{width=\figwidth}
    \includegraphics[width=\figwidth]{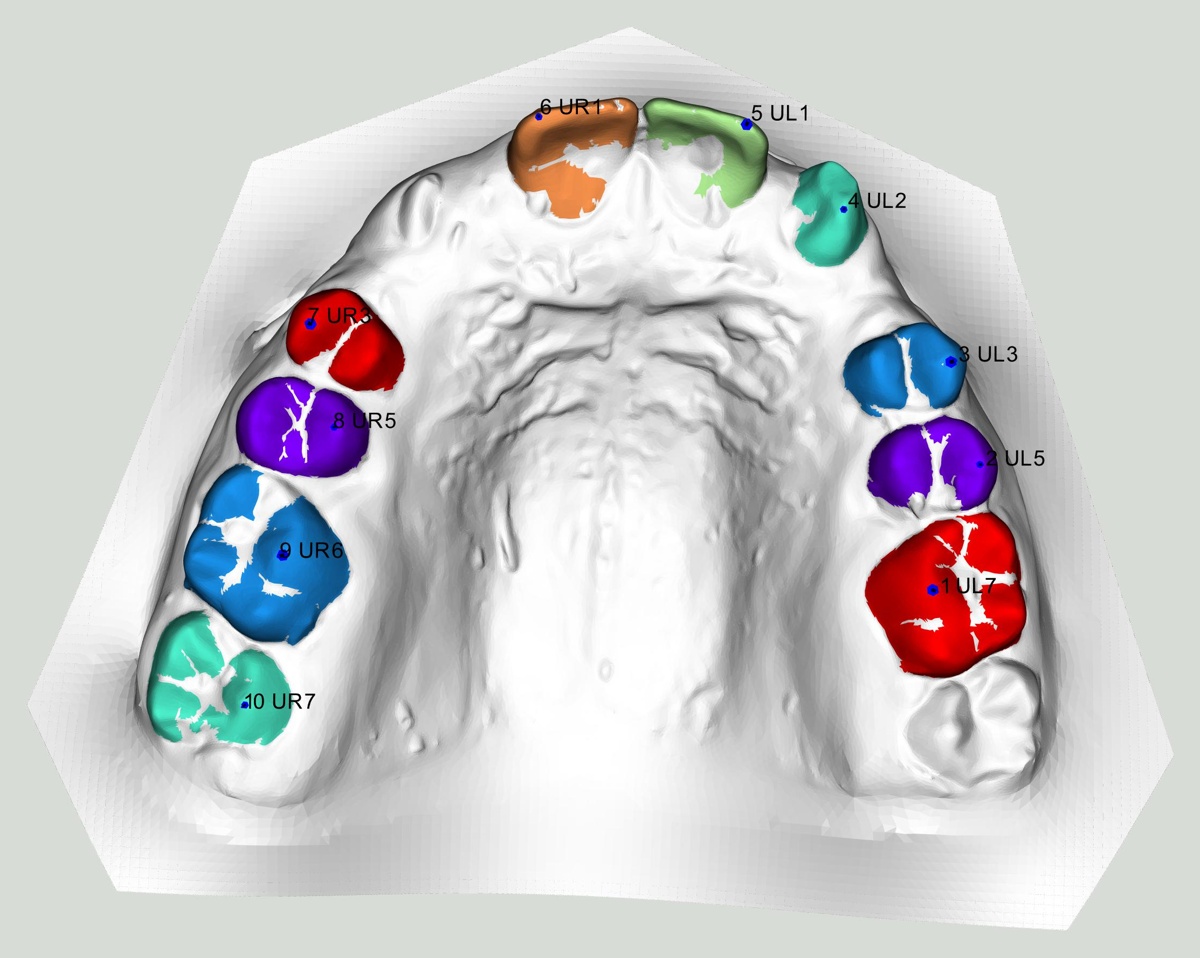}
    \caption{A model with several teeth misassigned. Section~\ref{sec:create-cost-table} needs some work.}
    \label{fig:assignment-guessing}
\end{figure}

The \nth{1} half of the assignment step (section~\ref{sec:create-cost-table}) is quite weak. For models with most, if not all, teeth present the linear programming (section~\ref{sec:solve-cost-table}) picks up the slack to give a good end result but for models with many teeth missing, such as the one shown in figure~\ref{fig:assignment-guessing}, the software gets progressively less reliable.

\begin{figure}[ht]
    \centering
    \captionsetup{width=\figwidth}
    \includegraphics[width=\figwidth]{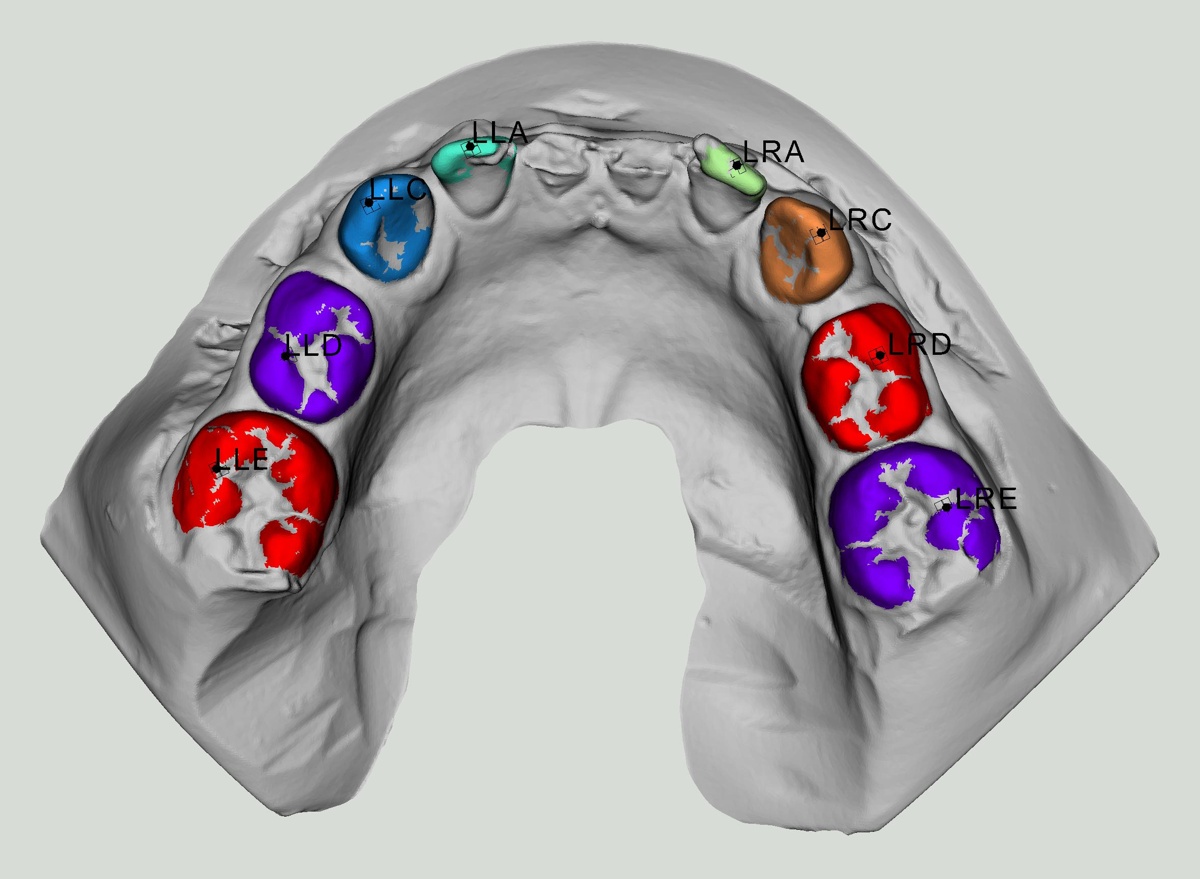}
    \caption{Neither the software, nor the person who programmed it, can tell lower incisors apart.}
    \label{fig:lower-incisors}
\end{figure}

The assignment works much better for deciduous models (albeit partly because it has less teeth to decide on) despite there being considerably more variation in deciduous teeth. Lower incisors are an exception (see figure~\ref{fig:lower-incisors}) but, given that even experienced clinicians struggle to distinguish lateral from central lower incisors, this is hardly surprising.

\begin{figure}
    \centering
    \captionsetup{width=\figwidth}
    \includegraphics[width=\figwidth]{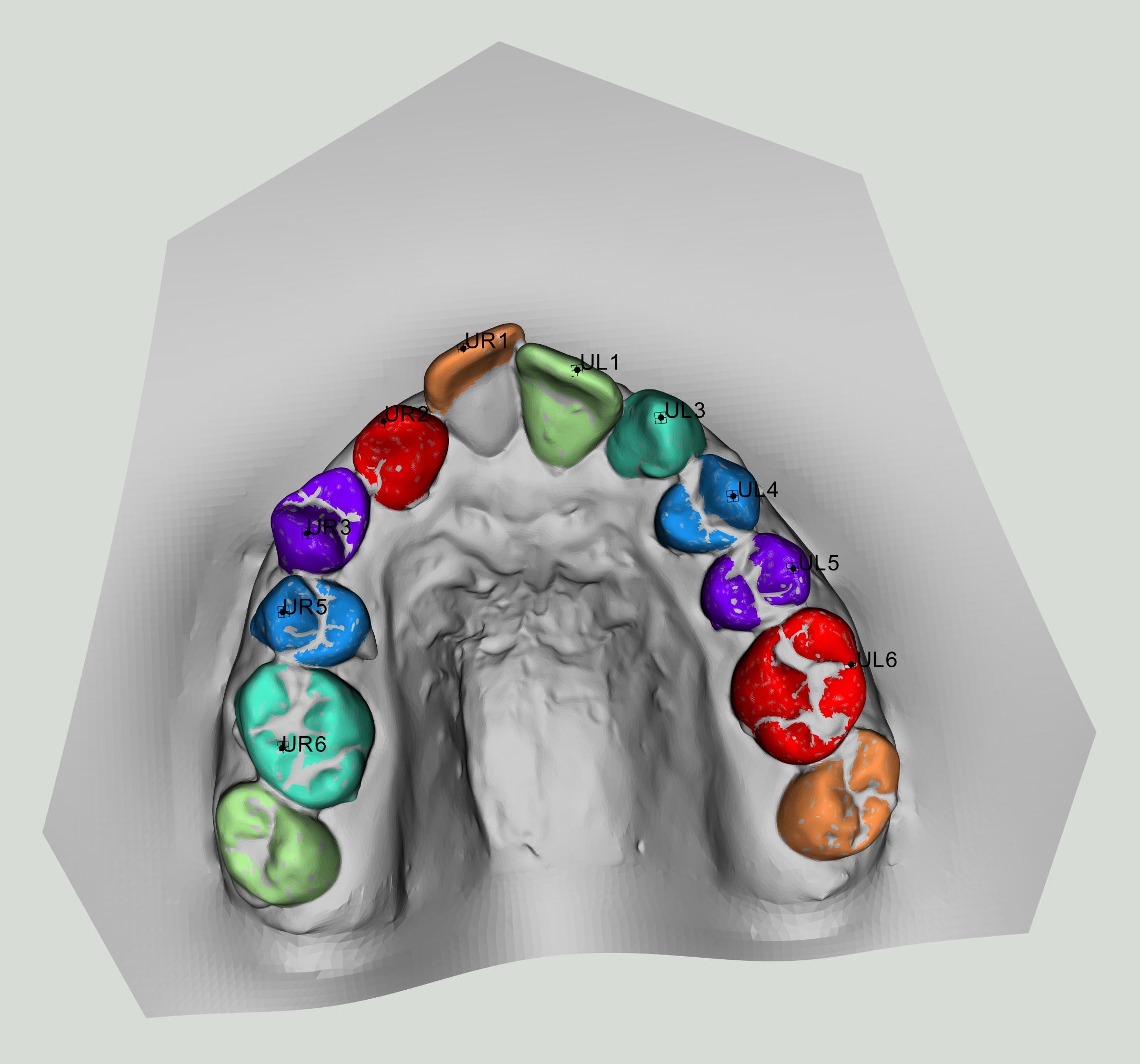}
    \caption{The orientation of the UR4 causes it to be mistaken for an UR3.}
    \label{fig:assignment-orientation-dependent}
\end{figure}

All but one of the tooth characteristics in section~\ref{sec:tooth-characteristics} are dependent on the orientation of the teeth. It's rare that a tooth be oriented unusually enough to cause a different assignment but it can happen as it did in figure~\ref{fig:assignment-orientation-dependent}. This dependence also rules out any chance of recognising teeth \emph{lose}, i.e. not part of a full model, which would be of great value to forensics. Ideally, characteristics should only measure orientation-independent properties such as curvature or normalise orientation first using something like PCA. 

\begin{figure}
    \centering
    \captionsetup{width=\figwidth}
    \includegraphics[width=\figwidth]{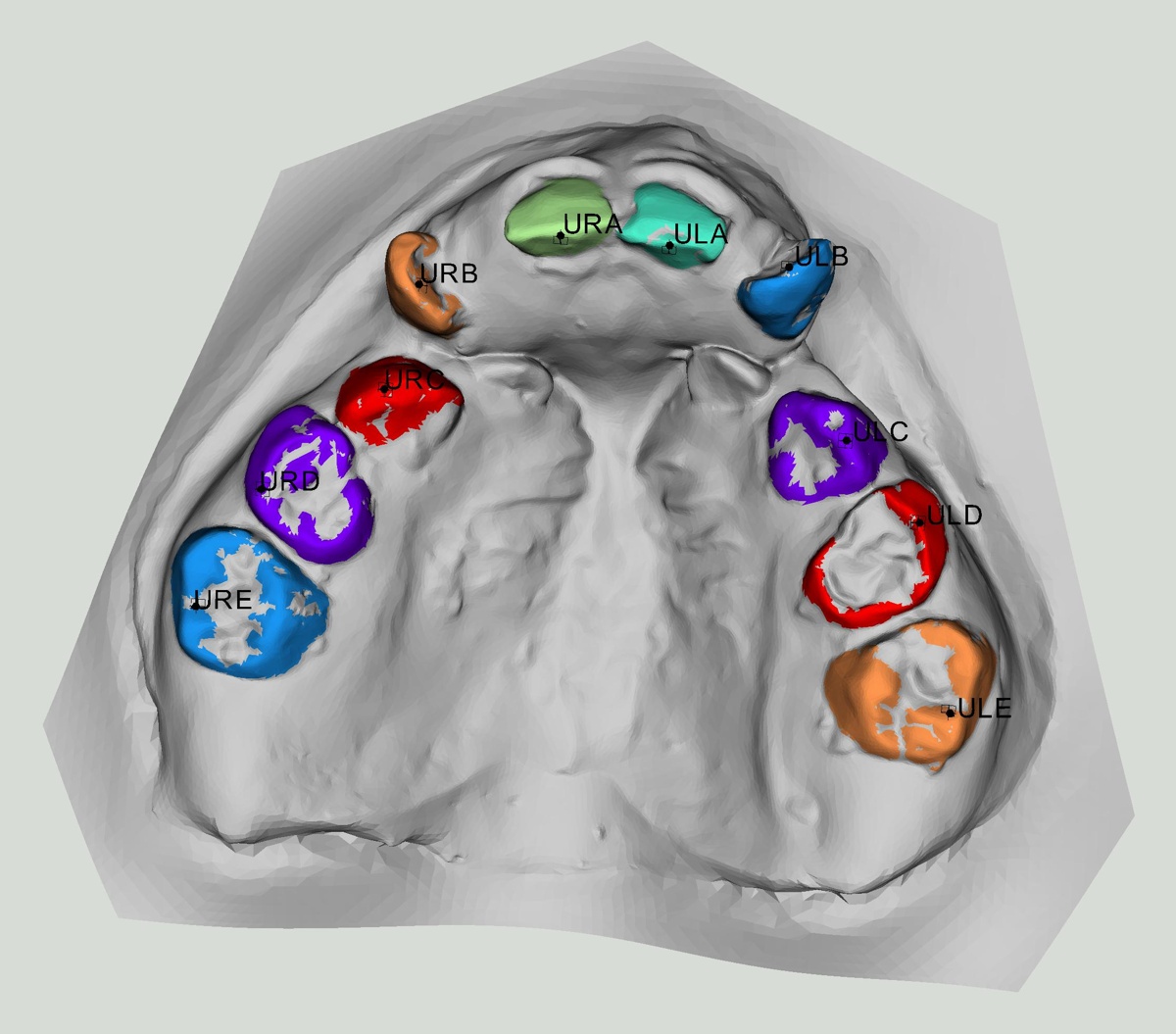}
    \caption{Despite the orientation or the two lateral incisors, linear programming is able to make the correct tooth assignments, albeit because it has very little choice.}
    \label{fig:assignment-orientation-dependent-addendum}
\end{figure}

With all teeth present, however, linear programming is again able to keep the final conclusion correct (see figure~\ref{fig:assignment-orientation-dependent-addendum}).

\section{Discussion}

Regarding tooth partitioning, one of the previous works showed an advantage in relation to ours. \cite{snake} used teeth outlines to achieve teeth segmentation. One of the complications that arised from our methodology, was that molar teeth were often recognised as two teeth instead of one entity. Whilst our software is later able to correct this issue after identifying each tooth, the segmentation procedure by \cite{snake} doesn't make this error to begin with. Our software doesn't truly partition teeth, often leaving gaps in grooves on molars or ignoring the lingual surface of incisors. For finding landmarks this is not an issue but for forensic dentistry, true partitioning is requirement to achieve accurate identification. True partitioning can be achieved by passing annotated landmarks from our software as control points to \cite{harmonic-field-ZOU2015132} if required.
\textit{Harmonic fields} \cite{harmonic-field-ZOU2015132} remains the gold standard of tooth partitioning with the drawback that it is only semi-automatic. Given that our own software, whilst technically fully automatic, still requires proof checking by hand, our software will not supersede until substantial improvement in reliability is achieved.

We had the opportunity to test our software on a considerably more varied dataset than any prior works we found. Our dataset has been crucial to ensuring our software will not perform drastically worse on all models par those used for development. Our models trickled in batches of 5-20 models. These model batches are data shared with us from other studies (of patients consented to have their data also used by us). As a result, each batch was very different from the last. The data sets included a wide range of patients' ages and conditions and two different brands of scanner. We also used plaster models from different centers that were created and scanned by many different clinicians. This is something I see other studies would have benefited from. If \cite{Automatic-Feature-Identification-in-Dental-Meshes} had also had tried scans from a TRIOS scanner, which uses the Y-axis for vertical, they would have known not to simply hard-code vertical as the Z-axis. Similarly, if \cite{snake} had had access to near-toothless models, certain intra-oral scans or plaster models with knobbly/textured bases, they would have discovered that raw PCA is distorted by such models and is insufficient to orientate with exclusively.

Similarly to how dentists require a single shareable means to benchmark treatments (a primary goal of this project), a common set of shared dentitions would be a requirement to benchmark the softwares of different researching groups. The success of tooth partitioning in particular greatly depends on the models you give it. The \textit{crispness} of edges, the smoothness of tooth surfaces, the presence of pockmarks and welts on plaster models or the use of intra-oral scanning all affect the quality of results. For all those reasons, it is very hard to meaningfully cross compare the reliability of our software techniques to others like it.


\section{Conclusion}
\label{sec:conclusion}

In this paper we develop from scratch a set of methods for the automated landmarks recognition from scanned mesh data of dental surfaces. The automatically identified landmarks are crucial for developing an automated scoring software based on the MHB system to measure the outcome of dental treatments. The specific requirements of the MHB system and the difficult to predict effects of complex geometry of patient teeth request original thoughts and innovative methods which are not readily available in literature. 

Our methods include the following steps: 
\begin{enumerate}
    \item Use the center of mass, principal component analysis and fitting a gradient line to tips, to find an approximate position and orientation of the created dental surface from its scanned data;
    \item Use the local maxima in the vertical direction to automatically provide an initial approximation of the landmarks; 
    \item Extract surface gradient and curvature information to identify the shape and boundaries of individual tooth -- developing a 3D image segmentation technique specific for the purpose of tooth segmentation; 
    \item Order teeth through a best-fit quadratic jaw-line approximation; 
    \item Use a combination of machine learning and linear programming to recognize and label each tooth and its landmarks. 
\end{enumerate}  

We also provide quite a few prior attempts that were tried but didn't work so as to prevent any future developers from making the same mistakes. We have successfully automated the MHB scoring system by using the methods studied in this paper. (The MHB software's details and scoring performance will be reported elsewhere in future.) Furthermore, this software has a much broader application. It can be expanded to automatically identify landmarks for a range of other scoring indices.


\bibliographystyle{abbrv}
\bibliography{references.bib}

\appendix

\section{Calculating PCA}
\label{sec:calculate-pca}

Information describing how to calculate PCA is rather sparsely available. Hence, a recipe to apply PCA to a set of points is included below.

Consider all points in the model.
\begin{equation*}
    X = \lbrace \mathbf{x_1}, \mathbf{x_2}, \dots, \mathbf{x_n} \rbrace
\end{equation*}
Subtract the centre of mass from each point to get displacements. The centre of mass being the mean of all points.
\begin{equation*}
    X' = \lbrace \mathbf{x_i} - \Bar{\mathbf{X}} \:, \; i=1 \dots n \rbrace
\end{equation*}

Matrix multiply $X'$ transposed with itself.
\begin{equation*}
    M = X'^T X'
\end{equation*}

$M$ should be a $3x3$ matrix and is typically referred to as the \emph{covariance matrix}.

Then use eigen decomposition on $M$ to get three eigenvalues and their corresponding eigenvectors. The eigenvectors are the unit-vectors / directions and the eigenvalues are the covariances in those directions. The unit-vectors should be sorted by eigenvalue from largest to smallest.

\end{document}